%% file: arxiv_main.tex
\pdfoutput=1
\documentclass[ba,preprint]{imsart}

\RequirePackage{amsthm,amsmath,amsfonts,amssymb}
\RequirePackage[authoryear]{natbib}
\RequirePackage[colorlinks,citecolor=blue,urlcolor=blue]{hyperref}
\RequirePackage{graphicx}

\pubyear{2026}
\arxiv{2603.04199}
\volume{TBA}
\issue{TBA}
\firstpage{1}
\lastpage{1}

\startlocaldefs
\theoremstyle{plain}

\newtheorem{theorem}{Theorem}[section]

\theoremstyle{definition}

\newtheorem{remark}[theorem]{Remark}

\theoremstyle{remark}

\input{preamble_v5}
\endlocaldefs

\begin{document}

\begin{frontmatter}
\title{Bayesian Adversarial Privacy\bapsubtitle{A Bayesian decision-theoretic framework for contextual privacy analysis}} 
\runtitle{Bayesian Adversarial Privacy}

\begin{aug}
\author[A]{\fnms{Cameron}~\snm{Bell}\ead[label=e1]{bell@ceremade.dauphine.fr}},
\author[A]{\fnms{Timothy}~\snm{Johnston}\ead[label=e2]{johnston@ceremade.dauphine.fr}},
\author[A]{\fnms{Antoine}~\snm{Luciano}\ead[label=e3]{luciano@ceremade.dauphine.fr}},
\and
\author[A,B]{\fnms{Christian P}~\snm{Robert}\ead[label=e4]{xian@ceremade.dauphine.fr}}


\address[A]{CEREMADE, Université Paris Dauphine--PSL, France%
\printead[presep={,\ }]{e1,e2,e3}}

\address[B]{Department of Statistics, University of Warwick, UK%
\printead[presep={,\ }]{e4}}

\runauthor{Bell, Johnston, Luciano and Robert}

\end{aug}

\begin{abstract}
\input{sections/abstract}
\end{abstract}




\end{frontmatter}

\section{Introduction}
\label{sec:intro}
\input{sections/introduction_v6.tex}

\section{Bayesian Adversarial Privacy}
\label{sec:framework}
\input{sections/agents_v5}

\section{Example 1: A Coin Toss}
\label{sec:coin-toss}
\input{sections/coin_toss_v6}

\section{Example 2: Gaussian Hypothesis Testing}
\label{sec:gaussian}
\input{sections/gaussian_example_v5}

\section{Conclusion}
\label{sec:conclusion}
\input{sections/conclusion}

\begin{acks}[Acknowledgments]
The authors are grateful to the participants in the second and third OCEAN privacy workshops in Les Houches, February 2025 and Venice, March 2026, for their input and comments. They also warmly thank the entire editorial team for their constructive comments and suggestions.
\end{acks}

\begin{funding}

Timothy Johnston and Christian P Robert are respectively fully and partly supported by the European Union under the ERC Synergy Grant 101071601 (OCEAN, 2023–2030). Views and opinions expressed are solely those of the authors and do not necessarily reflect those of the European Union or the European Research Council Executive Agency. Neither the European Union nor the granting authority can be held responsible for them. Cameron Bell and Antoine Luciano are fully supported by a PR[AI]RIE-PSAI Chair funded by the Agence Nationale de la Recherche (ANR-23-IACL-0008). Part of the revision was written while Christian P Robert was a visiting professor at the Ca' Foscari University of Venice.

\end{funding}

\begin{supplement}
\stitle{Supplementary Material}
\sdescription{The supplementary material of the journal version (integrated risk derivations for the Gaussian example) is reproduced below as an appendix, Sections~S1--S11. The code used to generate the figures and numerical results is available at \url{https://github.com/AntoineLuciano/BayesianAdversarialPrivacy}.}
\end{supplement}

\clearpage
\bibliographystyle{ba2}
\bibliography{Bibliography}

\clearpage
\setcounter{section}{0}\renewcommand{\thesection}{S\arabic{section}}
\setcounter{figure}{0}\renewcommand{\thefigure}{S\arabic{figure}}
\setcounter{table}{0}\renewcommand{\thetable}{S\arabic{table}}
\setcounter{equation}{0}\renewcommand{\theequation}{S\arabic{equation}}
\renewcommand{\theHsection}{S\arabic{section}}\renewcommand{\theHfigure}{S\arabic{figure}}\renewcommand{\theHtable}{S\arabic{table}}\renewcommand{\theHequation}{S\arabic{equation}}
\section*{Supplementary Material: Integrated Risk Derivations for the Gaussian Example}
\input{supp_body}

\end{document}

%% file: preamble_v5.tex
\usepackage[svgnames]{xcolor}
\usepackage{enumitem}   
\usepackage{centernot}
\usepackage{dsfont}
\newcommand{\mathbbm}{\mathds}
\usepackage{comment}

\DeclareRobustCommand{\bapsubtitle}[1]{\texorpdfstring{\\ {\small #1}}{: #1}}

\DeclareMathOperator*{\argmin}{arg\,min}

\DeclareMathOperator{\Var}{Var}
\DeclareMathOperator{\E}{\mathbb{E}}

\DeclareMathOperator{\Prob}{\mathbb{P}}
\DeclareMathOperator{\ind}{\mathds{1}}

\newcommand{\dif}{\mathrm{d}}
\newcommand{\N}{\mathcal{N}}

\usepackage{amsmath,amssymb}
\usepackage{tikz}
\usetikzlibrary{arrows.meta,positioning,fit,calc,shapes.geometric,backgrounds}

%% file: sections/abstract.tex
Theoretical and applied research into privacy encompasses an incredibly broad swathe of differing approaches, emphases and aims. This work introduces a novel quantitative notion of privacy that is both contextual and specific. Building on and extending ideas from statistical disclosure control and differential privacy, our aim is to model the implications of a disclosure decision in an adversarial setting. Our definition relies on concepts inherent to standard Bayesian decision theory, while departing from them in several important respects. In particular, (i) inference about the data itself becomes meaningful and (ii) the party controlling the release of sensitive information should make disclosure decisions from the prior viewpoint, rather than conditional on the data, which is a feature shared with Bayesian design. Illuminating toy examples are exploited towards highlighting the specificities of the method.




%% file: sections/introduction_v6.tex

    The notion of privacy is increasingly appearing as  a central issue in fields concerned with the gathering and handling of personal and sensitive data. Its explosion in the past decades has been followed by a progressive realisation of the dangers posed by independent entities given access to one's personal details. This concern naturally includes statistics, which must balance the conflicting constraints of statistical efficiency and data protection. While the notions of \emph{data protection} and \emph{right to privacy} are widely discussed, the majority of scholars involved in the domain adopt the position that privacy occurs within a rich web of convention and context \citep[see e.g.][]{Rachels_1984,Nissenbaum2004PrivacyAC,cofone2023privacy}, hence allowing for highly varying definitions and measurements. For instance, Differential Privacy (DP) \citep{dwork2006differential}---reviewed in Section \ref{sec:others}---is often considered as the privacy canon, but it is associated with a specific understanding of data protection, regardless of contextual interpretations and statistical objectives. Whenever statistical efficiency is concerned, it must manoeuvre within predetermined DP constraints that may prove orthogonal to statistical imperatives \citep{duchy:2013,butucea:issartel:2021}
    all while resting on parameters that are notoriously hard to interpret \citep{cunmings:2021}.


    When incorporating privacy protection within a genuine inferential framework, especially with a Bayesian perspective, the goals of this protection should be spelled out into a quantitative definition that is layered onto the traditional statistical model, rather than remaining an abstract, context-free notion detached from the statistical problem. As repeatedly demonstrated in the Bayesian literature, a natural, rational, and conclusive approach to this aim is to proceed within a decision-theoretic framework \citep{raiffa:1968,berger:1985,robert2001bayesian}, because spelling out utility functions exposes the ultimate goals of the analysis in quite an explicit manner. The construction we develop in this paper proceeds from this imperative: it models as utilities the most essential implications of disclosing information about observed data in a statistical analysis, that is, balancing the inferential benefits of disclosure against the privacy harms it induces.


    While Bayesian versions of DP have been proposed in the literature, the corresponding papers mostly build on alternative definitions \citep{Desfontaines2019SoKDP} or establish privacy guarantees for broad classes of algorithms \citep{chourasia2021advances}, rather than studying the trade-off between privacy constraints and statistical accuracy \citep{eigner2014differentially}. As illustrated in \cite{Dimitrakakis2017DifferentialPF}, the assumptions DP imposes on the likelihood function and on the underlying Markov Chain Monte Carlo (MCMC) sequence generally impede statistical efficiency, for instance by limiting the number of MCMC samples. By looking at worst-case scenarios and imposing randomisation and insufficiency, DP's guarantees rarely coincide with Bayesian optimality criteria. The same holds for Statistical Disclosure Control (SDC), also reviewed in Section \ref{sec:others}, which is practically verified but also 
    often formulated independently of an explicit Bayesian model, even if its risk assessment can be cast in Bayesian terms.
    
    In contrast, while the approach adopted in this paper may be deemed too formal in its requiring a utility function to be spelled out, it focuses on statistical objectives and 
    incorporates data as a new component of statistical inference, where the data play a dual role: an observation to condition upon and, simultaneously, an inferential target to be protected. 
    The Bayesian analysis that comes together with this framework is furthermore quite appealing in that it incorporates prior information on both the statistical parameter and the data to be protected through a Bayesian model. As in non-Bayesian statistical approaches set in DP \citep{duchy:2013}, our approach proves almost always unable to achieve full optimality in the classical Bayesian decision-theoretic sense \citep{berger:1985}. In both settings, theirs and ours, the difficulty stands with the statistical procedures being selected prior to observing the data, in contrast with Bayesian optimal design \citep{duarte:wong:2016}.
    Differential privacy is designed to \emph{certify} a mechanism as acceptable rather than to \emph{rank} or \emph{compare} the mechanisms that qualify: $\epsilon$-DP is a threshold many mechanisms meet, and distinguishing among them is not its purpose.
    However, Bayesian adversarial privacy allows for a complete comparison of privacy procedures (or {\em release mechanisms}) constructed from first principles, both deterministic and randomised. 
    It further opens a Bayesian perspective on inferring about the data to be protected, given the released inference. Indeed, this feature delivers a quantitative statement about the privacy degree of this release, both a priori when assessing releases and a posteriori when addressing data holders' concerns.



The goal of this paper therefore differs somewhat from those of the established canons. Relative to differential privacy, our formalism gives inference an explicit role and takes into account the consequences of adversarial actions on the released summaries, where DP is set independently of any inferential objective. Relative to SDC, the difference is smaller and perhaps better seen as a change of perspective: SDC constrains disclosure risk to an acceptable level and uses an information-loss measure to select among admissible releases, while our criterion instead attaches a price to that constraint through a single trade-off parameter. In neither case is our proposal antagonistic: it can be used to compare and rank mechanisms designed under either canon, and, as we show in Section~\ref{subsec:dp-comparison}, it connects to both, DP reappearing as a bound on the privacy risk and SDC as the constrained form of our criterion.

In this paper, Section \ref{sec:others} reviews the existing perspectives on privacy that relate to our goal. Section~\ref{sec:framework} introduces the full decision-theoretic framework for an actor who observes the data and designs a release mechanism balancing information efficiency and privacy protection. Sections \ref{sec:coin-toss} and \ref{sec:gaussian} illustrate this novel framework in simple models where computations can be carried out the entire way. Section \ref{sec:conclusion} elaborates on the extensions of the current work and the many challenges faced by realistic settings.

\section{Concurrent Perspectives on Statistical Privacy}\label{sec:others}

\subsection{Differential Privacy}
Much recent literature on data privacy has focused on DP \citep{dwork2006differential}, a perturbative framework that constrains the sensitivity of the released output distribution as a function of the private input data. Specifically, one supposes that the output is a randomised mechanism (e.g., noise may be added to the data), and that small changes in the data induce only small changes in the \emph{probability distribution} of the output. There has been a proliferation of variant definitions featuring slightly different specifications of this core idea \citep{Desfontaines2019SoKDP}. The original definition in \citet{dwork2006differential} states that a mechanism $q(\cdot\mid x)$ is $\epsilon$-DP if, for all datasets $x,x'$ differing in one entry, and all measurable sets $A$
\begin{equation}\label{eq: DP defn}
\frac{q(A \mid x)}{q(A \mid x')}
\;\le\;
e^{\epsilon}\, .
\end{equation}
This condition ensures that the inclusion or removal of any single individual has a limited effect on the (distribution of) the output. See \citet{balle2020hypothesis} for an explanation via the language of hypothesis testing.

DP definitions compose coherently under natural probabilistic operations \citep{dwork2014algorithmic}, which makes the notion a good base for studying trade-offs between algorithmic and statistical accuracy and privacy constraints \citep{chourasia2021advances, eigner2014differentially}, even though it does not intrinsically pertain to the latter. While differential privacy implies strong privacy guarantees on membership inference against individuals, by construction, it is indeed set independently of any prior information or statistical consideration. Furthermore, if a stronger statistical model is assumed on the data, DP guarantees can end up safeguarding privacy less than advertised \citep{kifer2011no, ghosh2016inferential}. In addition DP cannot protect against the use of information learnt from populations against individuals \citep{1273551}, which can represent a serious privacy threat.

Since DP relies on worst-case analysis, realistic cases may require adding a large quantity of noise, which negatively impacts statistical inference \citep{doi:10.1126/sciadv.abk3283}.

DP is now deployed at scale, most prominently in the 2020 US Census \citep{desfontaines2023list,Bailie2026Refreshment}, and these deployments show that the parameter $\epsilon$ of Equation~\eqref{eq: DP defn} accumulates under composition, with budgets that are large and hard to interpret \citep{drechsler2023government,drechsler2024complexities}; see \citet{del2025differential} for a more critical perspective.

As discussed in Section \ref{sec:intro}, connections with Bayesian statistics have been given in works such as \citet{Dimitrakakis2017DifferentialPF}, where posterior sampling is utilised as a random mechanism. Other Bayesian DP variants \citep{ghosh2016inferential, Kasiviswanathan2008OnT, Yang2015BayesianDP,hu2026contaminated} require a prior distribution on the data, as we do, and impose for the posterior distribution given the released statistic to remain close enough to the full posterior. A particular branch of the DP literature opposes privacy vs statistical accuracy as an (economic) trade-off \citep{GHOSH2015334, 10.1145/3139457}. This formulation bears strong similarities to our framework, though correlations between the data and the privacy protection requested for it may complicate the setup.
This difficulty has much in common with our selecting the release mechanism \emph{ex ante}, that is, integrated over the prior predictive distribution. On the other hand, since our priors and losses are outlined explicitly from the start, much of the discussion on choosing a pricing scheme is made redundant.

\subsection{Statistical Disclosure Control}

The alternative notion of SDC has been mainly developed by National Statistical Institutes (NSIs) to process official data to be published without leaking confidential information \citep{hundepool2024handbook}. The central objective of SDC is to balance \emph{disclosure risk} and \emph{data utility}, the latter typically quantified through an \emph{information-loss} measure. For instance, the marginal totals of tables and other summary statistics should be preserved, whilst simultaneously preventing record matching with other non-anonymised released information. For example, with only two companies from a given sector, publishing their averaged revenue enables each company to learn the revenue data of its competitor. 

The notions and methods of SDC have expanded in such a diffuse way amongst disparate teams of researchers that it is hard to succinctly review the state of the art. However, it should be noted that the emphasis is generally on ensuring that the published data is close to the true data in a chosen metric, whilst at the same time preventing realistic data attacks, as opposed to building a statistical model for the data and the release. To ensure stronger privacy guarantees, SDC often hides the exact nature of the disclosure method \citep{slavkovic2023statistical}. This can be a drawback for statisticians, since adding noise leads to inconsistencies between published tables and to an uncertainty that users can hardly quantify \citep{eu2024}. Furthermore, one can argue that if an attacker manages to infer the disclosure method they could learn more than intended. 

Although DP and SDC represent two different traditions for privacy analysis, a substantial literature connects them. For instance, \citet{Bailie2026Refreshment} considers the foundations of DP in a manner closer to SDC, specifying precisely what characteristics are protected and how. This fruitful compromise reveals the extent to which DP is not a catch-all formulation for privacy \citep{Schneider2025}. It should be noted that DP-style analysis can be far more pessimistic than SDC risk measures about a release mechanism \citep{2423656.2423658}, since DP is agnostic to the actual values in the data whereas probabilistic disclosure risks depend on them.
Section~\ref{subsec:dp-comparison} makes both connections precise: an $\epsilon$-DP guarantee bounds Eve's prior-averaged risk from below whatever the prior, and the SDC practice of maximising utility under a cap on disclosure risk is the constrained form of our criterion.

\subsection{Further Related Literature}
Bayesian assessments of statistical disclosure that relate to our goals are found in previous works. For instance, in \citet{Fienberg1997ABA} adversarial attempts to infer additional information about individuals from a randomised data release are assessed via simulations, in \citet{Hu2018BayesianEO}, methods of reconstructing real data from synthetic data are discussed, \citet{neunhoeffer:etal:2026} establishes the Bayesian equivalence between a Gaussian randomisation mechanism and a Gaussian synthetic data generator and in \citet{Reiter2014BayesianEO}, posterior odds of an adversary matching datasets are considered, albeit with a worst-case analysis. In addition, \citet{Eilat2021BayesianP} considers mechanisms with a restriction on the Kullback–Leibler divergence between the prior and posterior. 

Furthermore, statistical information has already been set in an adversarial or minimax sense. The concept of ``informational size'' relates to the extent agents' beliefs being modulated by information releases (see \citet{McLean2001InformationalSA} and related works), and ``rational inattention'' \citep{Sims2003ImplicationsOR} considers the problem of constructing a random variable $Y$ that yields information about another variable $X$ subject to an information flow restraint.

More broadly, federated learning, where data in separate silos is processed in such a way as to remain confidential, is a popular and active topic of research \citep{kairouz2021advances}. Extensive connections have been made with both Bayesian analysis \citep{cao2023bayesian} and DP \citep{wei2020federated}, and it has proved extremely effective in practical settings, such as with healthcare data \citep{Sheller2020}.

Quantitative information flow (QIF) \citep{alvim2020science} offers many connections with our approach. The idea of QIF is to model what is leaked by an information release (a \emph{channel}), in a Bayesian setting where the data is given a prior distribution. Many subtleties of this formulation are explored in detail in \cite{alvim2020science}, for instance sensitivity to the choice of prior, the impact on the choice of adversarial loss ($g$-\emph{vulnerability}), worst case scenario analysis over families of priors and adversarial losses, the composition of releases and the impact of hidden correlations (\emph{Dalenius leakage}).

Finally, both our proposal and the recent proposal by \citet{bon2026persuasiveprivacy} develop a closely related Bayesian decision-theoretic notion of privacy. In fact, both papers share the same representation of statistical data privacy towards aggregating the utility of conducting statistical inference with the loss resulting from privacy leakages. The common underlying agent-based, game-theoretic, framework opens a range of Bayesian privacy strategies. In the approach of \citet{bon2026persuasiveprivacy}, a scoring rule is chosen on the parameter---as opposed to data---predictive distribution and a supplementary constraint is imposed on the privacy level. More precisely, one considers a mechanism $q(\cdot\mid x)$ to be private enough if, for every predictive distribution associated with a family of priors and for every realisation $x$ in the data-space, the probability of facing a high score difference in an optimal inference between prior and posterior predictives stays low enough. This alternative perspective therefore adopts a more minimax, less Bayesian, version of privacy than ours, but remains connected with our decision-theoretic, multi-agent, framework, with the ensuing complexity that a worst-case solution need be found over a family of priors. However, one of the differentiating properties of their framework is that it recovers the original $(\epsilon, \delta)$-DP definition as a special case. In addition, several structural features of DP carry over to the wider notion, including post-processing and composition rules, as well as handling deterministic release mechanisms.

%% file: sections/agents_v5.tex
We consider a Bayesian setting in which, given data $x$, inference is to be performed on a parameter $\theta$ driving the distribution of $x$. Performing such inference requires releasing information derived from $x$, which may contain sensitive content. Our objective is to evaluate release mechanisms according to both the quality of inference on $\theta$ and the privacy leakage regarding $x$.

Let $\theta \in \Theta$ denote an unknown parameter indexing the distribution of the data and $x \in \mathcal{X}$ the observed data. We assume a prior probability measure $\pi(\dif \theta)$ on $\Theta$ and a likelihood kernel $p(\dif x \mid \theta)$ on $\mathcal{X}$ as a function of $\theta$. The joint distribution driving the Bayesian analysis is therefore given by $\pi(\dif \theta, \dif x) = p(\dif x \mid \theta)\,\pi(\dif \theta)$, with marginal distribution $p(\dif x) = \int_\Theta p(\dif x \mid \theta)\,\pi(\dif \theta)$, and posterior distribution $\pi(\dif \theta \mid x) = \pi(\dif\theta) p(x \mid \theta)/p(x)$. We assume that all distributions have a density with respect to some reference measures. Note that due to the specifics of Bayesian adversarial privacy, $p(\cdot)$ will also be called prior predictive (or sometimes just prior) distribution on $x$. An extension of our framework to a prior distribution on $x$ differing from $p(\cdot)$ for privacy attack purposes can easily be constructed.

Although we focus on parameter inference, the framework applies equally when the data themselves are the primary object of interest for both parties, namely for scenarios such as census data. In that case, one may integrate over $\theta$ to obtain a marginal formulation.

We will be considering the prior, likelihood and loss functions as fixed throughout. However, it is of course of great practical interest to investigate the behaviour of different release mechanisms under a variety of models or losses. As in any Bayesian setting, these choices are modelling assumptions and should be considered to reflect the interests of the user, rather than informing some underlying absolute truth. See \cite{bon2026persuasiveprivacy} for an alternative approach that operates over a family of prior distributions.

\subsection{Three Agents with Linked Goals}

To formalise and expose the trade-off between inference and privacy, we introduce three imaginary agents to represent the competing objectives, Alice, Bob and Eve, named after the famous couple of the RSA foundational paper \citep{RSA:78}, as well as an early privacy paper \citep{BBR:1988}.

\textbf{Alice} is the \emph{mechanism designer} (e.g. a server). She selects a release mechanism $q$, generating a release $\eta \in \mathcal{H}$ according to $\eta \sim q(\dif \eta \mid x)$, depending on the data $x$. The release induces a posterior distribution on both the data and the parameter, which we write as
\[
p(\dif x \mid \eta, q)
=
\frac{q(\eta \mid x)\,p(\dif x)}
{\int_{\mathcal{X}} q(\eta \mid x')\,p(\dif x')}, 
\qquad 
\pi(\dif \theta \mid \eta,q) = \int_{\mathcal{X}} \pi(\dif \theta \mid x)\,p(\dif x \mid \eta,q),
\]
to explicitly include the dependence on $q$. 

The chosen mechanism $q$ and the observed $\eta$ are assumed to be available to the public, while $x$ remains known only to Alice. Alice wishes to perform inference on $\theta$ without compromising characteristics of $x$ she wants to hide. To do this, she will have to choose $q$ in such a way that observing $\eta$ yields useful evidence for $\theta$, whilst protecting elements of interest in $x$. To this end, she brings in two hypothetical agents, Bob and Eve, who will represent Alice's dual objectives.

\textbf{Bob}, the \emph{statistician}, aims to infer the parameter $\theta$ using only the released output $\eta$ and knowledge of the release mechanism $q$. Let $L_B(\theta,\delta)$ denote his loss function with decision space $\mathcal{D}_B$. Given $(\eta,q)$, Bob's posterior expected loss and optimal Bayes decision are then
\begin{equation}
\rho_B((\eta,q),\delta)
=
\int_\Theta L_B(\theta,\delta)\,\pi(\dif \theta \mid \eta,q), \qquad
\delta_B(\eta,q)
=
\arg\min_{\delta \in \mathcal{D}_B}
\rho_B((\eta,q),\delta).
\label{equ-delta-Bob}
\end{equation}

The choice of $L_B$ determines what inference is to be done, for example point estimation, hypothesis testing, prediction, or uncertainty quantification. For instance, quadratic loss leads to posterior mean estimation, whereas $0$–$1$ loss yields posterior mode decisions in testing problems \citep{robert2001bayesian}. For Alice, $L_B$ could represent her own inferential goals (making Bob truly hypothetical), or could correspond to those of a real client.

Meanwhile, \textbf{Eve}, the \emph{adversary} or \emph{eavesdropper}, seeks to learn about the data $x$ from $(\eta,q)$. With loss $L_E(x,\delta)$ and decision space $\mathcal{D}_E$, her posterior expected loss and Bayes rule are
\begin{equation}
\rho_E((\eta,q),\delta)
=
\int_{\mathcal{X}} L_E(x,\delta)\,p(\dif x \mid \eta,q),
\qquad
\delta_E(\eta,q)
=
\arg\min_{\delta \in \mathcal{D}_E}
\rho_E((\eta,q),\delta).
\label{equ-delta-Eve}
\end{equation}

Similarly to $L_B$, the specification of $L_E$ determines which features of the dataset Alice views as sensitive.  The loss function formalises the notion of a privacy breach, whether this concerns the presence of outliers, the value of specific individual observations, or broader structural properties of the data. Unlike Bob, it is clear that Eve will always be a hypothetical construct: a genuinely malicious actor would not preemptively alert a server of what they were after. For this reason, it is of utmost importance for Alice to understand exactly what features of $x$ she wishes to keep private. This is a big difference with DP, for example, where privacy of a mechanism is a much more universal property with respect to the data.

\begin{remark}[Multiple Bobs or Eves]
The specificity of $L_B$ and $L_E$ raises the following natural question: what should Alice do if there are ``multiple Bobs or Eves'', each representing a different inference or privacy goal of Alice's? Such a setting can be placed into the current framework as follows. Focusing on ``multiple Eves'' (as the same method can be used for multiple Bobs), suppose each Eve has a loss function $L_E^{(i)}$, for $i= 1, \dots, n$, each of which has associated Bayes rule $\delta_E^{(i)}$. Alice can balance these different adversaries by linearly combining them into a single loss function
\begin{equation*}
    L_E(x, (\delta^{(1)}, \dots, \delta^{(n)}) ) = \sum_{i=1}^n w_i L_E^{(i)}(x, \delta^{(i)}),
\end{equation*}
with positive weights $w_i$ for each loss. It is clear that $\delta_E = (\delta_E^{(1)}, \dots, \delta_E^{(n)})$ is an optimal vector decision for this loss, which may not hold if the losses were combined using, say, a minimum. The relative weights $\{ w_i\}_{i=1}^n$ would need to be carefully specified by the user, factoring in potential heterogeneity of scale between the various losses. In this sense, multiple Bobs or Eves are really equivalent to a single agent with a more complicated objective.
\end{remark}

\subsection{Alice's Objective}
\label{section-alice}

Within this framework, Alice needs to select the release mechanism $q$. Her task is to balance two competing goals: improving Bob’s inference on $\theta$ while limiting Eve's inference on $x$.

We define Alice’s loss as the expected difference between Bob’s and Eve’s losses, evaluated at their respective Bayes-optimal decisions under the mechanism $q$:
\begin{equation*}
L_A((\theta,x), q)
=
\int_{\mathcal H}
L_B\!\left(\theta, \delta_B(\eta,q)\right)
\, q(\dif \eta \mid x)
-
\lambda
\int_{\mathcal H}
L_E\!\left(x, \delta_E(\eta,q)\right)
\, q(\dif \eta \mid x),
\end{equation*}
where $\lambda>0$ determines the relative weight placed on statistical utility and privacy protection. Alice's loss can be seen as a version of Bob's loss with a ``privacy penalty'' on the release, similar to a Lasso or ridge penalty. She aims to use as much of $x$ to improve Bob's inference on $\theta$ as possible, without releasing aspects of the data that would allow Eve to make a better decision $\delta_E$ about $x$.


In the special case where the release is deterministic, so that $\eta=\eta(x)$, Alice’s loss reduces to
\begin{align*}
L_A((\theta,x), q)
=
L_B\!\left(\theta,\delta_B(\eta(x),q)\right)
-
\lambda
L_E\!\left(x,\delta_E(\eta(x),q)\right).
\end{align*}

To evaluate the quality of a release mechanism $q$, Alice should (for reasons we will highlight below) consider the Bayesian \emph{integrated mechanism risk}, also known as \emph{ex ante risk}, given by
\begin{equation}
\begin{split}
R_A(\pi,q)
&=
\iint
L_A((\theta,x),q)
\,\pi(\dif \theta,\dif x)
\\[6pt]
&=
\iint
L_B(\theta,\delta_B(\eta,q))
\,q(\dif \eta \mid \theta)
\,\pi(\dif \theta)
\\
&\quad
-
\lambda
\iint
L_E(x,\delta_E(\eta,q))
\,q(\dif \eta \mid x)
\,p(\dif x).
\end{split}
\label{eq:alice_risk}
\end{equation}
This ex ante risk is widely used in Bayesian experimental design and related fields \citep{Chaloner1995BED,FriedExAnte2003,Rainforth2024}. This criterion averages Alice’s loss over the joint model and corresponds to the Bayes risk adapted to our adversarial setting, providing a principled and coherent basis for comparing mechanisms. One could attempt to find a global minimiser of $R_A(\pi,q)$ over all $q$'s, but in general this will be intractable, as we shall see below.

A crucial observation is that $L_A((\theta,x),q)$ depends not only on the conditional distribution $q(\cdot \mid x)$, but on the entire mechanism $q(\cdot \mid y)$ for all $y \in \mathcal X$. In other words, Alice must consider not just what her output is at a specific observed data $x$, but also what her outputs could be for every possible $y$. This dependence arises through the induced posterior distribution $p(\dif x \mid \eta,q)$, since the decision rules $\delta_B$ and $\delta_E$ minimise integrals against this density. In particular, two mechanisms $q_1$ and $q_2$ may induce identical distributions of the release at the realised dataset while yielding different losses, that is to say $q_1(\cdot \mid x) \stackrel{D}{=} q_2(\cdot \mid x)$ does not imply $L_A((\theta,x),q_1) = L_A((\theta,x),q_2)$. 

To illustrate this distinction, let $\mathcal X=\{0,1\}$ and consider two deterministic mechanisms: a full release $\eta_1(x)=x$ and a constant $\eta_2(x)=0$. Then $q_1(\cdot \mid x=0)=q_2(\cdot \mid x=0)$, so the releases coincide at $x=0$. However, the induced posterior distributions differ:
\[
p(\dif x \mid \eta=0,q_1)=\delta_0(\cdot),
\qquad \text{whereas} \qquad
p(\dif x \mid \eta=0,q_2)\neq \delta_0(\cdot).
\]
Since Bob and Eve optimise against these posterior distributions, the resulting decisions $\delta_B$ and $\delta_E$, and hence Alice’s loss, generally differ. For this reason, a seemingly protective release such as “release $x$ unless it would help Eve, in which case release nothing” may fail to provide privacy. Since Eve knows the mechanism $q$, observing “nothing” itself conveys information, namely, that the realised dataset lies in a region deemed sensitive.

The global nature of $L_A$ distinguishes our problem from classical Bayesian decision theory. One might attempt to determine an optimal release locally at the observed dataset $x$ via
\begin{equation*}
    q^\star(\cdot \mid x)
    =
    \argmin_q
    \left(
        \int_\Theta
        L_A((\theta,x),q)
        \,\pi(\dif \theta \mid x)
    \right).
\end{equation*}
However, this approach is inherently circular. Even if $q^\star(\cdot \mid y)$ is fixed for all $y \neq x$, modifying $q^\star(\cdot \mid x)$ changes the induced posterior distributions $p(\dif x \mid \eta,q^\star)$, and therefore alters $L_A((\theta,y),q^\star)$ elsewhere in the sample space. Thus, optimality cannot be established pointwise, and we must rely on the ex ante risk.

Instead, Alice will typically consider a family of releases, possibly indexed by a parameter. Notable choices for $q$ include:
\begin{itemize}
    \item[--] The \emph{full release} where Alice releases the full dataset, so that $\eta = x$ a.s., and Bob has access to the full posterior $\pi(\dif \theta \mid x)$ for $\theta$. Meanwhile, Eve knows $x$ and her integrated risk reduces to the minimal value of $\min_\delta L_E(x,\delta)$.
    \item[--] The \emph{null release} where Alice releases nothing, which we denote by $\eta = \dagger$, so that Bob and Eve both rely exclusively on their respective prior information. We use $\dagger$ to refer to a cemetery state containing no information.
    \item[--] A \emph{random sample} from the full posterior $\pi(\dif \theta \mid x)$ for $\theta$, for instance produced by MCMC.
    \item[--] A \emph{noisy release}, for example $\eta = x + \epsilon$, for some random $\epsilon$.
    \item[--] A \emph{statistic} $\eta = T(x)$, such as the sample mean $\bar{x}$ or the posterior probability of an event in $\theta$.
    \item[--] A \emph{synthetic} dataset generated from the posterior predictive.
    \item[--] Different releases on different parts of $\mathcal{X}$, for example $\eta = \bar{x} + \epsilon$ when $x \in A$, and $\eta = \dagger$ when $x \notin A$.
\end{itemize}

\begin{remark}[Approximating $R_A$ via Monte Carlo]
\label{rmk:mc}
Evaluating $R_A$ could require costly, high-dimensional integrals, which may require the use of Monte Carlo methods to evaluate. Simulating tuples $(\theta, x, \eta)$ should be simple, but solving the optimisation problems for $\delta_B$ and $\delta_E$ from Equations~\eqref{equ-delta-Bob} and \eqref{equ-delta-Eve} at each $\eta$ could be costly, as these require first approximating the relevant posterior ($\pi(\dif\theta\mid\eta,q)$ for Bob, $p(\dif x\mid\eta,q)$ for Eve). In full generality, this can be done crudely via Approximate Bayesian Computation (ABC) and stochastic gradient methods, reusing the original $(\theta,x,\eta)$ samples to reduce the computational burden. In practice this fully nested optimisation is rarely needed, since standard losses give Bayes rules and risks that are explicit posterior functionals, e.g.\ for squared error losses that return posterior moments.
\end{remark}

Returning to our three hypothetical agents, Alice’s mechanism risk admits the decomposition
\begin{equation}
\label{eq:RA_decomposition}
R_A(\pi,q)
=
R_B(\pi,q)
-
\lambda R_E(\pi,q),
\end{equation}
where
\begin{equation*}
R_B(\pi,q) 
=
\iint
L_B(\theta,\delta_B(\eta,q))
\, q(\dif\eta\mid\theta)\,
\pi(\dif\theta),
\end{equation*}
and
\begin{equation*}
R_E(\pi,q)
=
\iint
L_E(x,\delta_E(\eta,q))
\, q(\dif\eta\mid x)\,
p(\dif x) .
\end{equation*}
The term $R_B(\pi,q)$, which we call the \emph{integrated inference risk}, is the expected loss induced by Bob’s optimal decision, averaged jointly over the prior uncertainty on $\theta$ and the randomness of the release $\eta$. Similarly, the term $R_E(\pi,q)$, the \emph{integrated privacy risk}, measures the expected loss of Eve’s optimal decision, averaged over the prior predictive distribution of the data $x$ and the released output $\eta$.

\begin{remark}[Misaligned agents]
\label{rmk:misaligned}
Bob and Eve are extensions of Alice rather than autonomous parties: they act under her prior, and $L_B$ and $L_E$ encode exactly the objectives she attributes to them. One may nonetheless ask what happens if the actual adversary departs from Eve, by holding a different prior $p_E\neq p$ or by pursuing a different objective $\tilde L_E\neq L_E$. In either case the adversary ends up playing some rule $\tilde\delta:\mathcal H\to\mathcal D_E$, whose harm Alice measures by her own standards through

\[
  R_E(\pi,q;\tilde\delta)=\iint L_E\bigl(x,\tilde\delta(\eta)\bigr)\,q(\dif\eta\mid x)\,p(\dif x).
 \]
Since $\delta_E(\eta,q)$ minimises $\rho_E((\eta,q),\cdot)$ at every $\eta$ by Equation~\eqref{equ-delta-Eve}, $R_E(\pi,q;\tilde\delta)\ge R_E(\pi,q)$ for every $\tilde\delta$: as measured by $L_E$, no adversary acting on $(\eta,q)$ is more harmful than Eve, and $R_E(\pi,q)$ is a worst-case privacy assessment. 

However, when moving to comparing mechanisms, the tightness of this bound may change the ranking between them. Whenever Alice can attribute a specific prior or objective to that adversary, she is therefore better advised to encode it in $p_E$ or $L_E$, at the price of one more posterior to compute, than to rely on the bound.
\end{remark}

\subsection{The Privacy-Utility Trade-off}
\label{subsec:ru-map}

The two integrated risks $R_B(\pi,q)$ and $R_E(\pi,q)$ admit a geometric representation in the \emph{risk--utility (R--U) confidentiality map} of \citet{duncan2004disclosure}. Writing $(U,R)=(-R_B(\pi,q),-R_E(\pi,q))$, each mechanism $q$ corresponds to a point whose first coordinate increases with the quality of Bob's inference and whose second coordinate increases with the success of Eve's attack. The full release then occupies the upper-right extreme of the attainable region and the null release its lower-left extreme (Figure~\ref{fig:ru-map}). On their own, these coordinates describe a mechanism but do not order mechanisms; they yield a decision criterion for Alice only once combined into the single objective $R_A=R_B-\lambda R_E$ of Equation~\eqref{eq:RA_decomposition}. Among the possible ways of combining them, this linear one is particularly convenient, as it reduces the trade-off to choosing the single parameter $\lambda$, to which we now turn.

This weight $\lambda>0$ governs the strength of the privacy penalty and encodes the normative balance between statistical utility and privacy protection. At $\lambda=0$, the criterion reduces to $R_B(\pi,q)$ and is minimised by the full release, inference being pursued with no regard for disclosure. As $\lambda\to\infty$, privacy dominates and any mechanism that conveys usable information to Eve is eventually beaten by the null release. Between these extremes, $\lambda$ fixes the rate at which Alice is willing to concede inferential utility for privacy protection. Its meaning is relative to the scales of $L_B$ and $L_E$: if the two losses are not commensurate, the induced trade-off is distorted, so that any meaningful comparison of mechanisms presupposes an explicit choice of $\lambda$.

A natural default calibration of Equation~\eqref{eq:RA_decomposition} is obtained by choosing $\lambda$ such that the full release $q^{\mathrm{full}}$ and the null release $q^\dagger$ yield identical integrated risk. Solving
\(
R_A(\pi,q^{\mathrm{full}})
=
R_A(\pi,q^\dagger)
\)
gives
\begin{equation}
\lambda
=
\frac{R_B(\pi,q^\dagger) - R_B(\pi,q^{\mathrm{full}})}
{R_E(\pi,q^\dagger) - R_E(\pi,q^{\mathrm{full}})},
\label{equ-lambda-full=dagger}
\end{equation}
whenever the denominator is non-zero. This places the two extremes at indifference, so that every non-trivial release is assessed relative to both, and other values of $\lambda$ correspond to different priorities between utility and privacy. 
This calibration is not the only reading of $\lambda$: in Section~\ref{subsec:dp-comparison} the same multiplier reappears, from a dual standpoint, as the shadow price of privacy in a disclosure-control problem.

\input{sections/ru_map}

The calibration \eqref{equ-lambda-full=dagger} has a direct geometric reading. Since $R_A=R_B-\lambda R_E=-(U-\lambda R)$, the level sets $\{R_A=\mathrm{const}\}$ form a family of parallel lines of slope $1/\lambda$ in the $(U,R)$ plane, along which $R_A$ decreases towards the bottom-right. The multiplier thus acquires a concrete meaning as the common slope of Alice's indifference lines. Requiring $q^{\mathrm{full}}$ and $q^{\dagger}$ to share the same risk, as in Equation~\eqref{equ-lambda-full=dagger}, is precisely the choice that places both on a single such line; equivalently, $\lambda$ is the reciprocal slope of the chord joining the null and full releases. It is against this line that every other mechanism is read: Alice prefers releases lying below it and, under her criterion, rejects any lying above it, such as the point $q'$ of Figure~\ref{fig:ru-map}.

Two displacements away from the extremes are especially telling, each an axis-parallel move on the map. The first leaves Eve's risk at its null level: a release that gives Eve no usable information leaves her posterior, and hence her Bayes decision, unchanged, so that she remains at the null risk $R_E^{\dagger}$ whatever the value of $\eta$. The second holds Bob's risk at its full-release floor $R_B^{\mathrm{full}}$: a release that leads Bob to the very decision he would have taken from the complete data preserves his risk exactly. Performing both at once lands at the \emph{corner} $(U,R)=(-R_B^{\mathrm{full}},-R_E^{\dagger})$, which couples the utility of the full release with the privacy of the null one. Under the calibration \eqref{equ-lambda-full=dagger}, and provided the full release lets Eve reconstruct her target exactly, so that $R_E^{\mathrm{full}}=0$, Alice's risk at the corner takes the value 
\begin{equation}
    R_A^{\star}=2R_B^{\mathrm{full}}-R_B^{\dagger},
    \label{equ-corner-general}
\end{equation}
the lowest value of $R_A$ attainable by any mechanism. Eve is absent from this value only because $\lambda$ has already absorbed her null-to-full privacy loss, not because she plays no role. Whether the bound is \emph{attained} depends squarely on her: a single release must saturate both bounds at once, which the correlation between the two objectives generally prevents, so the corner is best read as a benchmark reached only in stylised settings.

The vertical displacement admits a purely statistical construction. What Alice needs, in order to hold Bob at his full-release floor, is not that the release preserve his entire posterior but merely that it preserve his \emph{decision}: any $q$ with $\delta_B(\eta,q)=\delta_B(x,q^{\mathrm{full}})$ leaves $R_B$ unchanged and moves straight down the map, a property we call \emph{Bob-invariance}. Sufficiency for $\theta$ is one way to secure it, and a familiar one, since if $T(x)$ is sufficient then $\pi(\dif\theta\mid x)=\pi(\dif\theta\mid T)$, so the release $\eta=T(x)$ leaves Bob's posterior, decision and risk untouched, $R_B(\pi,q^{T})=R_B^{\mathrm{full}}$, while Eve's posterior no longer collapses to a point mass at $x$, whence $R_E(\pi,q^{T})\ge R_E(\pi,q^{\mathrm{full}})$ and $R_A(\pi,q^{T})\le R_A(\pi,q^{\mathrm{full}})$. The release of a sufficient statistic therefore weakly dominates the full release, whatever Eve's target and loss; sufficiency is nonetheless more than is required, as it preserves the whole posterior rather than only the decision it induces. 

The minimal Bob-invariant release is the decision itself, $\eta=\delta_B(x,q^{\mathrm{full}})$: conditioning on $\eta=d$ restricts to $\{x:\delta_B(x,q^{\mathrm{full}})=d\}$, on which $d$ is already optimal for every $x$, hence optimal on average.
Being coarser than any sufficient statistic, it can only help Eve less, and it is this release that reaches the corner whenever publishing $\delta_B(x)$ leaves Eve at her prior risk $R_E^{\dagger}$; the extent of the vertical descent is otherwise governed by how much of $x$ the released summary still betrays, as measured by $L_E$. This singles out a release that should always belong to Alice's admissible range, Bob's own decision, and connects with the use by \citet{fearnhead:prangle:2012} of the posterior mean as an optimal summary statistic in semi-automatic ABC. Example~1 realises the horizontal displacement (Section~\ref{sec:coin-toss}) and Example~2 the corner itself (Section~\ref{sec:gaussian}).

\subsection{Recovering Differential Privacy and Statistical Disclosure Control}
\label{subsec:dp-comparison}

Both canons reviewed in Section~\ref{sec:others} connect to our framework, though differently: differential privacy re-emerges as a \emph{bound} on Eve's risk under strong assumptions, while statistical disclosure control coincides with the constrained form of Alice's criterion.

Consider first differential privacy. Eve is an adversary distinguishing two neighbouring datasets $x\in\{D,D'\}$ under a uniform prior $p(D)=p(D')=\tfrac12$ and a $0$--$1$ loss $L_E(x,\delta)=\ind\{\delta\neq x\}$. If the mechanism $q$ is $\epsilon$-DP, with $e^{-\epsilon}q(\eta\mid D')\le q(\eta\mid D)\le e^{\epsilon}q(\eta\mid D')$, Bayes' rule gives $$\tfrac{1}{1+e^{\epsilon}}\le P(D\mid\eta)\le\tfrac{1}{1+e^{-\epsilon}},$$ so that Eve's integrated privacy risk satisfies $R_E(\pi,q)=\E[\min\{P(D\mid\eta),P(D'\mid\eta)\}]\ge\tfrac{1}{1+e^{\epsilon}}$, which coincides with the hypothesis-testing interpretation of DP \citep{balle2020hypothesis}. DP therefore imposes a \emph{uniform} lower bound on $R_E$ across every neighbouring pair and every prior, whereas BAP reports the actual prior-averaged $R_E$ of which the DP bound is merely the worst-case envelope. 
In the language of quantitative information flow, $-R_E(\pi,q)$ is the posterior $g$-vulnerability of the release for the gain $g=-L_E$, and this bound is a special case of the bounds that $\epsilon$-DP places on $g$-leakage, whatever the prior and the gain \citep[Chapter~23]{alvim2020science}.

Formulations from statistical disclosure control are recovered not as a bound but as a constrained program, and the connection runs through the multiplier $\lambda$ itself. In SDC, it is common to cap disclosure risk at an acceptable level and, under that cap, find a release minimising the loss of statistical utility, which is typically quantified through an information-loss measure. Writing the utility loss as $R_B$ and the cap as a floor $R_E\ge\varepsilon$, this is the program 
\[ \min_{q}\ R_B(\pi,q) \qquad\text{subject to}\qquad R_E(\pi,q)\ge\varepsilon, \] 
whose Lagrangian differs from Alice's criterion only by a constant. 
Minimising $R_A=R_B-\lambda R_E$ is therefore equivalent to solving this SDC-type problem, provided the attainable region of mechanisms is convex, so that a supporting multiplier exists, with $\lambda$ attached to the disclosure-risk constraint \citep{bertsekas:1982}. The multiplier then carries the interpretation of a \emph{shadow price of privacy}: the marginal inferential loss conceded for one additional unit of confidentiality, equivalently the slope $\dif R_B/\dif R_E$ of the risk--utility frontier of Section~\ref{subsec:ru-map}, whose reciprocal is the slope of the level lines of $R_A$ in Figure~\ref{fig:ru-map}.
What sets the Bayesian reading apart from standard SDC is that both $R_B$ and $R_E$ are \emph{ex ante} risks, expectations under the prior and the model rather than quantities evaluated at the realised data $x$ or the single release $\eta$: the disclosure-risk cap is thereby placed on a coherent probabilistic footing.

%% file: sections/ru_map.tex
\begin{figure}[t]
    \centering
    \definecolor{rumapSuff}{HTML}{6B87AC}
    \definecolor{rumapCorner}{HTML}{BF5B52}
    \definecolor{rumapBad}{HTML}{B07A2E}
    \definecolor{teal}{HTML}{1F8A8A}
    \begin{tikzpicture}[>=Stealth,line join=round,scale=0.95]
    \draw[->,gray] (0,0)--(8.8,0);
    \draw[->,gray] (0,0)--(0,6.8);
    \node[gray,below] at (4.4,-0.15){\small Utility $U=-R_B$};
    \node[gray,rotate=90,above] at (-0.2,3.4){\small Disclosure risk $R=-R_E$};

    \draw[densely dotted,gray] (6.7,5.75)--(6.7,1.0);
    \draw[densely dotted,gray] (1,1)--(6.7,1.0);

    \draw[teal!50,dashed,thick] (0.55,2.40)--(4.40,5.61);

    \draw[teal,very thick] (1,1)--(6.7,5.75)
        node[pos=0.48,sloped,above,font=\footnotesize]{%
            $R_A=R_A^{\dagger}=R_A^{\mathrm{full}}$\;(slope $1/\lambda$)};

    \draw[teal!50,dashed,thick] (3.70,0.20)--(7.70,3.53);

    \draw[->,black!70] (5.1,3.65)--(5.8,2.95);
    \node[black!70,right] at (5.5,3.5){\footnotesize $R_A\!\downarrow$};
    
    \draw[->,rumapSuff,thick,densely dotted] (1.20,1.0)--(3.75,1.0);
    \node[rumapSuff,below,align=center,font=\footnotesize] at (2.6,0.9){$R_E=R_E^{\dagger}$};

    \draw[->,rumapSuff,thick,densely dotted] (6.7,5.55)--(6.7,2.88);
    \node[rumapSuff,right,font=\footnotesize] at (6.86,4.70){$R_B=R_B^{\mathrm{full}}$};

    \fill[teal] (1,1) circle(2.8pt);
    \node[teal,anchor=north] at (1,0.55){\footnotesize Null $q^{\dagger}$};
    \fill[teal] (6.7,5.75) circle(2.8pt);
    \node[teal,anchor=south west] at (6.95,5.95){\footnotesize Full $q^{\mathrm{full}}$};
    \fill[rumapBad] (2.05,3.65) circle(2.8pt);
    \node[above left=1pt,rumapBad] at (2.05,3.65){\footnotesize $q'$};
    \fill[rumapSuff] (6.7,2.70) circle(3.0pt);
    \node[right,rumapSuff,font=\footnotesize] at (6.88,2.48){$q^{T}$};
    \draw[black!55,thick] (6.7,1.0) circle(2.7pt);
    \node[right,black!60,font=\footnotesize] at (6.9,0.6){corner $R_A^{\star}$};
    \fill[rumapSuff] (3.90,1.0) circle(3.0pt);
    \node[rumapSuff,above,font=\footnotesize] at (3.95,1.12){$q^{H}$};
\end{tikzpicture}
    \caption{A generic risk--utility map under the calibration of Equation~\eqref{equ-lambda-full=dagger}, each mechanism represented by the pair $(U,R)=(-R_B,-R_E)$. Dashed lines are level sets of $R_A$; the solid line is the level set through both the null and full releases. The point $q'$ above this line is ranked below both extreme releases. Below it, $q^{T}$ modifies the full release by publishing a statistic $T$ sufficient for $\theta$, which leaves Bob's risk unchanged (vertical displacement), while $q^{H}$ modifies the null release without informing Eve (horizontal displacement). The \emph{corner} marks the lowest attainable value $R_A^{\star}$ of Equation~\eqref{equ-corner-general}.}
    \label{fig:ru-map}
\end{figure}

%% file: sections/coin_toss_v6.tex
To illustrate the tension between statistical utility and privacy, and the importance of evaluating mechanisms ex ante, we consider a simple toy example: a medical study in which the recorded binary outcome is itself the sensitive condition, so that inference about its prevalence conflicts with protecting the individual record.

We model this with two coins, a fixed coin that always gives tails ($\theta=0$) and a fair coin ($\theta=1/2$), each chosen with prior probability $1/2$. Alice observes a single toss $X\mid\theta\sim\mathrm{Bern}(\theta)$, so that $X=1$ implies $\theta=1/2$. Bob infers the coin type $\theta$, Eve the realised toss $X$. They use decision spaces $\mathcal D_B=\{0,1/2\}$, $\mathcal D_E=\{0,1\}$ and losses
\[
L_B(\theta,\delta)=\ind\{\delta\neq\theta\},
\qquad
L_E(x,\delta)=
\begin{cases}
0, & \delta=x,\\
1, & \delta=1,\ x=0,\\
10, & \delta=0,\ x=1.
\end{cases}
\]
The asymmetry in $L_E$ makes a missed sensitive event ($x=1$) ten times costlier than a false positive, so Eve declares $\delta=0$ only when the posterior odds of $x=0$ exceed ten (the value $10$ is a modelling choice expressing that Alice deems the exposure of a sensitive record far costlier than a false alarm). We evaluate each mechanism through its integrated risk (Equation~\eqref{eq:alice_risk}), from the two trivial extremes to a globally optimal one.

\subsection{Full and Null Releases}

We start by considering the two trivial releases in order to bracket the trade-off and calibrate $\lambda$.

Under the \emph{full release} $q^{\mathrm{full}}$ ($\eta=x$), Eve recovers $X$ exactly, so $R_E=0$ and Alice's risk is independent of $\lambda$. Bob's posterior is $\Prob(\theta=0\mid X=0)=2/3$, so his Bayes rule is $\delta_B=1/2$ when $x=1$ and $\delta_B=0$ when $x=0$. He errs only when the fair coin lands tails and so mimics the fixed coin, that is when $\theta=1/2$ and $X=0$; hence
\[
R_A(\pi,q^{\mathrm{full}})=\Prob(\theta=1/2,\,X=0)=\tfrac14.
\]

Under the \emph{null release} $q^{\dagger}$ ($\eta=\dagger$), both agents act on the prior. Bob's two actions are then equally good (risk $1/2$), while Eve, by the asymmetry of her loss, declares $\delta_E=1$ and incurs $L_E=\ind\{x=0\}$; hence
\[
R_A(\pi,q^{\dagger})=\tfrac12-\tfrac34\lambda.
\]
The two extremes carry equal risk when $\lambda=\tfrac13$, which we adopt throughout.

\subsection{Randomised Response}

Between the two extremes lies a one-parameter family that randomises the released outcome, similarly to the approach of \cite{warner1965randomized}. For $\omega\in[0,1]$, the randomised response mechanism $q^{(\omega)}$ flips the outcome before releasing it, $\eta=X$ with probability $1-\omega$ and $\eta=1-X$ with probability $\omega$; both agents know $\omega$. Observing $\eta$ under $q^{(\omega)}$ is equivalent to observing $1-\eta$ under $q^{(1-\omega)}$, so $R_A$ is symmetric about $\omega=1/2$, and we may restrict without loss of generality to $\omega\in[0,1/2]$, where Alice tells the truth more often than not. On this range the family runs from the full release ($\omega=0$, always truthful) to the null release ($\omega=1/2$, pure noise), so its minimiser performs at least as well as either trivial release.

Both Bayes rules are explicit. Since Alice is more likely to tell the truth, Bob follows the signal, deciding $\theta=1/2$ exactly when $\eta=1$:
\begin{equation*}
    \delta_B(\eta,\omega)=\tfrac12\,\ind\{\eta=1\}.
    \label{equ-coin-toss-bob-decision}
\end{equation*}
Eve keeps her prior declaration $\delta_E=1$ and departs from it only when the signal is decisive, at $\eta=0$ once the flip rate falls below $3/13$ (a clean $\eta=0$, barely corrupted, points convincingly to $x=0$):
\begin{equation*}
    \delta_E(\eta,\omega)=
    \begin{cases}
        0, & \eta=0 \text{ and } \omega<\tfrac{3}{13},\\[4pt]
        1, & \text{otherwise}.
    \end{cases}
    \label{equ-coin-toss-eve-decision}
\end{equation*}
The threshold $3/13$ is where Eve's two posterior expected losses coincide; for $\omega\ge 3/13$ the signal is never decisive enough to move her, so her risk stays flat. This \emph{plateau}, which follows from the asymmetry of $L_E$, is what Alice exploits. Direct enumeration of the joint events $(\theta,X,\eta)$ gives
\begin{equation*}
R_B(\omega)=\tfrac14+\tfrac12\omega,
\qquad
R_E(\omega)=
\begin{cases}
\tfrac{13}{4}\,\omega, & 0\le\omega<\tfrac{3}{13},\\[6pt]
\tfrac34, & \tfrac{3}{13}\le\omega\le\tfrac12.
\end{cases}
\end{equation*}

Minimising $R_A(\pi,q^{(\omega)})=R_B(\omega)-\lambda R_E(\omega)$ over $\omega$ then gives
\begin{equation*}
\omega^\star(\lambda)=
\begin{cases}
0, & \lambda<\tfrac{2}{13},\\[6pt]
\tfrac{3}{13}, & \lambda\ge\tfrac{2}{13},
\end{cases}
\qquad
R_A(\pi,q^{(\omega^\star)})=
\begin{cases}
\tfrac14, & \lambda<\tfrac{2}{13},\\[6pt]
\tfrac{19}{52}-\tfrac34\lambda, & \lambda\ge\tfrac{2}{13}.
\end{cases}
\end{equation*}
At the calibrated $\lambda=1/3$, the optimum is $\omega^\star=3/13$, the boundary of Eve's plateau: the release carries strictly more information about $X$ than the null release, improving Bob's inference on $\theta$, yet the asymmetry of $L_E$ leaves Eve's Bayes action, and her risk $R_E=3/4$, unchanged. Randomised response thus exploits all the information irrelevant to Eve's decision before any of it benefits her. Even when privacy is heavily weighted, calibrated randomisation strictly dominates both extremes.

\begin{remark}[A differential-privacy reading]
The randomised response family is $\epsilon$-differentially private by construction. Since $q^{(\omega)}(\eta\mid x)$ only takes the values $1-\omega$ and $\omega$, the likelihood ratio between the two data values is bounded by $(1-\omega)/\omega$, so $q^{(\omega)}$ satisfies $\epsilon$-DP with
\[
\epsilon(\omega)=\log\frac{1-\omega}{\omega},
\]
decreasing from $\epsilon=\infty$ at the full release ($\omega=0$) to $\epsilon=0$ at the null release ($\omega=\tfrac12$). Minimising $R_A$ over the family therefore selects an operating point on the $\epsilon$-DP frontier: under $\lambda=\tfrac13$ the optimum $\omega^\star=\tfrac{3}{13}$ corresponds to $\epsilon^\star=\log\tfrac{10}{3}\approx1.20$. This illustrates, on a tractable example, how BAP can act as a selection principle over a family of differentially private mechanisms (Section~\ref{subsec:dp-comparison}).
\end{remark}

\subsection{Optimal Solution via Linear Programming}
\label{subsec:linear_prog}

The collection of randomised response is only a one-parameter subfamily. To optimise over all mechanisms, we exploit that the parameter, sample and decision spaces are finite, which recasts Alice's problem as a finite Linear Program (LP). From any release $\eta$ Alice can compute the Bayes decisions $\delta_B,\delta_E$ deterministically, so without loss of generality she releases the pair $\eta=(\delta_B,\delta_E)\in\mathcal D_B\times\mathcal D_E$, parametrised by the joint probabilities $\rho_{ij}^k=\Prob(\eta=(i,j),X=k)$. Writing $\tilde L_B(x,i)=\int L_B(\theta,i)\,\pi(\dif\theta\mid x)$ for the posterior-averaged inference loss, Alice's risk
\begin{equation*}
R_A(\pi,q^{(\rho)})
=
\sum_{k}\sum_{i}\sum_{j}
\Bigl(
\tilde L_B(k,i)-\lambda L_E(k,j)
\Bigr)
\rho_{ij}^k
\label{eq:RA_linear_rho}
\end{equation*}
is linear in $\rho$. The simplex constraints $\sum_{i,j}\rho_{ij}^k=\Prob(X=k)$, together with the requirement that $(i,j)$ be Bayes-optimal for both agents (a set of linear inequalities in $\rho$), define a finite-dimensional linear program, solvable exactly in the binary case.

Here the optimum is $q(\eta=(0,1)\mid X=0)=1$ and, for $X=1$,
\[
q(\eta=(0,1)\mid X=1)=\tfrac{3}{10},
\qquad
q(\eta=(1/2,1)\mid X=1)=\tfrac{7}{10}.
\]
The two signals play distinct roles. The signal $(1/2,1)$ occurs only under $X=1$, revealing both $X=1$ and $\theta=1/2$. The signal $(0,1)$ pools all $X=0$ outcomes with the remaining $30\%$ of the $X=1$ outcomes; conditional on it $\Prob(X=0\mid\eta)=\tfrac{10}{11}$, so Eve's two actions incur equal loss and she is exactly indifferent (playing $\delta_E=1$ only by tie-breaking), while $\theta=0$ remains Bob's posterior mode ($\Prob(\theta=0\mid\eta)=\tfrac{20}{33}$). Rather than actively misleading Eve, the mechanism places her posterior exactly on her decision boundary. The resulting risk
\[ 
R_A(\pi,q^{(\rho^\star)}) = 
\begin{cases} 
\tfrac{1}{4}, & \lambda \le \tfrac{1}{10},\\[6pt]
\tfrac{13}{40}-\tfrac{3}{4}\lambda, & \lambda > \tfrac{1}{10},
\end{cases} 
\] 
is strictly below the best randomised response ($R_B=\tfrac{13}{40}$ versus $\tfrac{19}{52}$, at the same $R_E=\tfrac34$). The improvement comes neither from adding noise nor from a richer signalling alphabet, since the optimum also uses only two signals; it comes from dropping the \emph{symmetry} of randomised response and instead tailoring the channel to the two agents' distinct losses. The mechanism is built so that the agents' Bayes decisions do exactly what Alice intends. Under $X=0$ she always sends $(0,1)$, which leads Bob to the correct decision; under $X=1$ she sends the truthful $(1/2,1)$ with probability $7/10$ and, with probability $3/10$, the pooled signal $(0,1)$ that steers Bob to the wrong decision $\delta_B=0$. Alice therefore misleads Bob only when $X=1$, and since $\Prob(X=1)=\tfrac14$ she does so on $\tfrac14\cdot\tfrac{3}{10}=\tfrac{3}{40}$ of datasets, exactly the $7.5\%$ by which $R_B=\tfrac{13}{40}$ exceeds the full-release floor $\tfrac14$. This is the price of pinning Eve's posterior on her indifference boundary: her action $\delta_E=1$ is the same under both signals, so she is never misled, merely kept uninformed.

\subsection{Summary of the Coin-Toss Example}

Fixing the calibration $\lambda=1/3$, for which $R_A(\pi,q^{\mathrm{full}})=R_A(\pi,q^{\dagger})=1/4$, Table~\ref{tab-coin-summary} reports the risks of each mechanism.

\begin{table}[ht]
    \centering
    \begin{tabular}{c|ccc}
        Mechanism & $R_B$ & $R_E$ & $R_A$\\
        \hline
        Full & 0.250 & 0.000 & 0.250 \\
        Null  & 0.500 & 0.750 & 0.250 \\
        Randomised ($\omega^\star = \tfrac{3}{13}$) & 0.365 & 0.750 & 0.115 \\
        Linear-program optimum & 0.325 & 0.750 & \textbf{0.075}
    \end{tabular}
    \caption{Integrated risks $R_B$, $R_E$ and $R_A$ for the coin-toss example with $\lambda=1/3$. Values are exact and reported to three decimals.}
    \label{tab-coin-summary}
\end{table}

\begin{figure}[ht]
    \centering
    \includegraphics[width=0.7\textwidth]{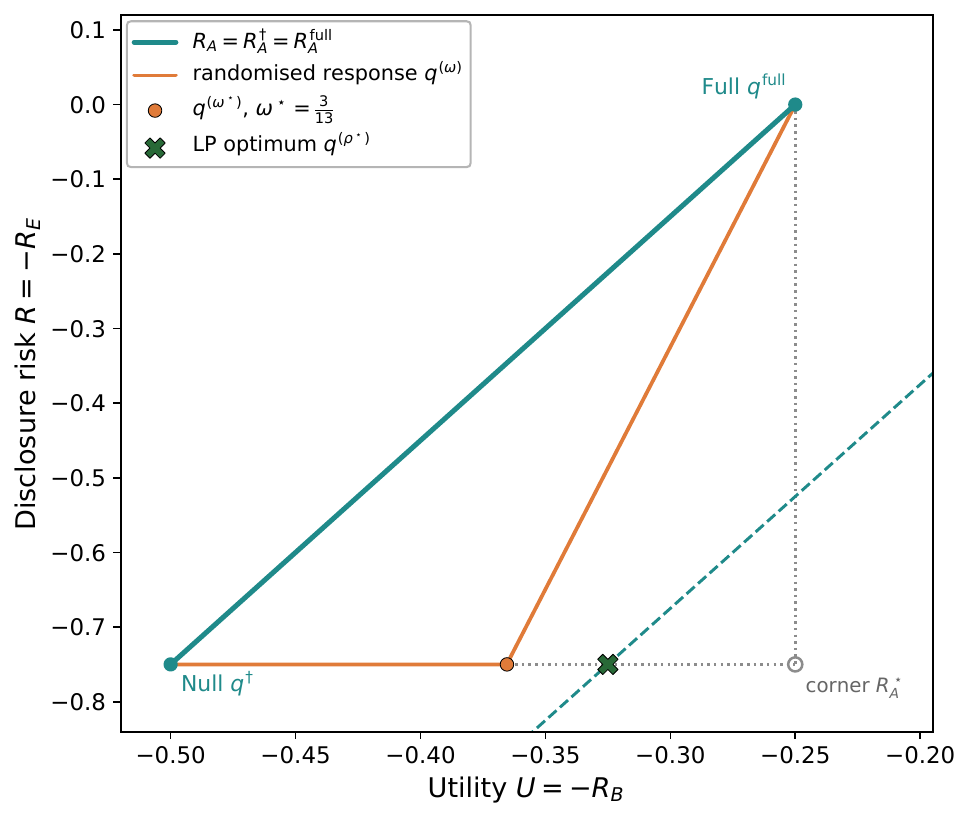}
    \caption{Risk--utility map of the coin-toss example: each mechanism is the point $(U,R)=(-R_B,-R_E)$, and the level set of $R_A$ through the full and null releases is the calibration line.}
\label{fig:rumap-cointoss}
\end{figure}

Figure~\ref{fig:rumap-cointoss} places these mechanisms in the risk--utility plane. Although the full and null releases are calibrated to share the same risk, both are strictly dominated. Randomised responses improve on them by exploiting the plateau in Eve's Bayes rule, moving away from non-disclosure only up to the boundary at which she would change her decision; the linear program improves further by abandoning symmetry, pooling the two data states asymmetrically so that one induced posterior sits exactly on Eve's indifference boundary, thereby preserving $R_E$ while conveying more decision-relevant information to Bob.

%% file: sections/gaussian_example_v5.tex
As a second example, we consider a Bayesian hypothesis testing problem in a Normal–Normal model. 
The aim is to investigate how inference on the mean parameter $\theta$ can be balanced against protecting different structural features of the dataset.

We adopt the conjugate Gaussian model
\begin{equation*}
    \begin{gathered}
        \theta \sim \N(0,\sigma_0^2),\\
        X=(X_1,\dots,X_n)\mid \theta \sim \N(\theta,1),
    \end{gathered}
\end{equation*}
and write $\bar X := n^{-1}\sum_{i=1}^n X_i$ for the sample mean.


Given thresholds $c_B, c_E \in \mathbb{R}$, we first assume Bob is interested in the null and alternative hypotheses
\begin{equation*}
    H_{0,B}:\ \theta \leq c_B, \qquad \mathrm{and} \qquad H_{1,B}:\ \theta > c_B,
\end{equation*}
respectively. Similarly, we assume Eve is interested in testing the hypotheses
\begin{equation*}
    H_{0,E} : \ T(X) \leq c_E, \qquad \mathrm{and} \qquad H_{1,E} : \ T(X) > c_E,
\end{equation*}
where we consider two choices for the statistic $T(X)$: $T(x) = \bar{x}$ in Section~\ref{sec-gaussian-mean}, or $T(x) = \max_i x_i$ in Section~\ref{sec-gaussian-max}. These two tests correspond to Eve wanting to identify either the centre of the sample or the presence of extreme values in the sample. 
The case where Eve targets both features at once is treated in the appendix (Section~S8).

We expect that the case where Eve is interested in $\bar{X}$ will be more challenging for Alice, since $\bar{X}$ is a sufficient statistic for $\theta$. Information about the mean will therefore be directly tied to the posterior in $\theta$ that Bob will be investigating. If Eve is interested in $\max_i X_i$, which is a statistic of the tails of the sample, Alice can potentially find stronger releases that help Bob without helping Eve.


Both agents use a $0$–$1$ loss,
\begin{equation*}
    L_B(\theta,\delta)
    =
    \ind\{\delta \neq \ind\{\theta>c_B\}\},
    \qquad
    L_E(x,\delta)
    =
    \ind\{\delta \neq \ind\{T(x)>c_E\}\}.
\end{equation*}
Hence, their optimal decisions are posterior majority rules, respectively
\begin{equation*}
    \delta_B(\eta,q)
    =
    \ind\{\pi(H_{1,B}\mid \eta,q)>1/2\},
    \qquad
    \delta_E(\eta,q)
    =
    \ind\{p(H_{1,E}\mid \eta,q)>1/2\}.
\end{equation*}


We compare the following mechanisms:
\begin{enumerate}[label=(\roman*)]

    \item \textit{Full release}: 
    $q^{\mathrm{full}}(\cdot\mid x)=\delta_x(\cdot)$ (so $\eta=x$).
    \vspace{0.4em}

    \item \textit{Null release}: 
    $q^\dagger(\cdot\mid x)=\delta_\dagger(\cdot)$ 
    (so $\eta=\dagger$ is independent of $x$).
    \vspace{0.4em}

    \item \textit{Noisy full release}: 
    $q^{\mathrm{full}}_\sigma$, releasing $Y=x+\varepsilon$ with
    $\varepsilon\sim \N(0,\sigma^2 I_n)$ independent, i.e.
    $Y_i\mid x_i\sim\N(x_i,\sigma^2)$.
    \vspace{0.4em}

    \item \textit{Noisy mean release}: 
    $q^{\mathrm{mean}}_\sigma$, releasing 
    $\eta=\bar x+\xi$ with 
    $\xi\sim \N(0,\sigma^2)$ independent.
    \vspace{0.4em}

    \item \textit{Noisy median release}: 
    $q^{\mathrm{med}}_\sigma$, releasing 
    $\eta=\mathrm{med}(x_{1:n})+\xi$ with 
    $\xi\sim \N(0,\sigma^2)$ independent.
    \vspace{0.4em}

    \item \textit{One-bit release}: 
    $q^{1\mathrm{bit}}_\tau$, releasing 
    $\eta=\mathbbm{1}\{\pi(H_{1,B}\mid x)>\tau\}$ 
    for $\tau\in[0,1]$.

\end{enumerate}

This is not an exhaustive list of possible mechanisms, but note that the noisy mean release is equivalent to releasing a posterior sample of $\theta$'s by sufficiency arguments. The same applies to the release of one synthetic dataset from the posterior predictive.

Additionally, releases (iii)-(vi) outlined above are parameterised releases, with the noisy releases depending on noise level (standard deviation) $\sigma$ and the one-bit release depending on a threshold $\tau$. We can therefore seek to minimise $R_A$ over the different releases, then compare these based on their minimised integrated risk.




The class of noisy full releases can be seen to contain both the full release (when $\sigma = 0$) and the null release (as the limit case when $\sigma \rightarrow \infty$). Both other calibrated releases have similar inference-privacy trade-offs, but do not necessarily contain the full release. Note furthermore that certain release calibrations recover Bob's optimal risk, for instance when $q$ is equal to the one-bit release with $\tau = 1/2$.

Throughout this section $\lambda$ is calibrated as in Equation~\eqref{equ-lambda-full=dagger}, so that the null and full releases share the same integrated risk. Therefore, following Equation~\eqref{equ-corner-general} one observes the general \emph{corner} of the risk--utility map is
\begin{equation}
    R_A^{\star}=2R_B^{\mathrm{full}}-\min(\pi(H_{0,B}),\pi(H_{1,B})),
    \label{equ-corner-test}
\end{equation}
which depends only on Bob's problem. However, whether this corner is actually \emph{attained} depends on the dependence between Bob and Eve's objectives. We exploit this in the design of the experiment by choosing the tests and the statistics they act on: below we fix Bob's threshold $c_B$ and Eve's two thresholds $c_E$ so as to place the mean adversary just short of the corner and the maximum adversary on it.

For the numerical experiments, we fix $n=5$, $\sigma_0=1$ and $c_B=0$. Eve's two thresholds sit one half and two prior standard deviations out, $c_E^{\mathrm{mean}}=\sigma_0/2$ and $c_E^{\max}=2\sigma_0$, respectively. 
The same two thresholds are reused for the multiple-Eve criterion of the appendix (Section~S8).
Explicit derivations of the posterior distributions and the corresponding optimal decisions 
are provided in 
the appendix.

\subsection{Privacy against the Mean}
\label{sec-gaussian-mean}
We first consider the case where Eve targets the sample mean, $T(x)=\bar x$. 

Since $\bar X$ is a sufficient statistic for $\theta$, Bob’s and Eve’s testing problems depend on essentially the same information. Consequently, improving Bob’s inference necessarily tends to reduce Eve’s uncertainty, leading to a genuine inference–privacy trade-off in the ex ante risk.

\begin{figure}[ht]
    \centering
    \begin{minipage}[b]{0.49\textwidth}\centering
        \includegraphics[width=\textwidth]{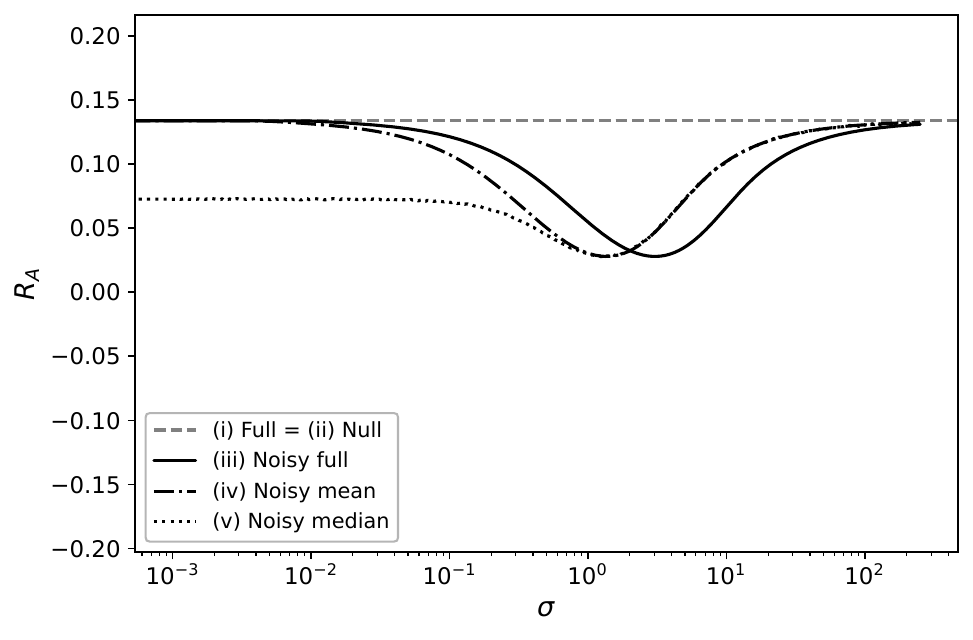}\\{\footnotesize (a) additive-noise releases vs.\ $\sigma$}
    \end{minipage}\hfill
    \begin{minipage}[b]{0.49\textwidth}\centering
        \includegraphics[width=\textwidth]{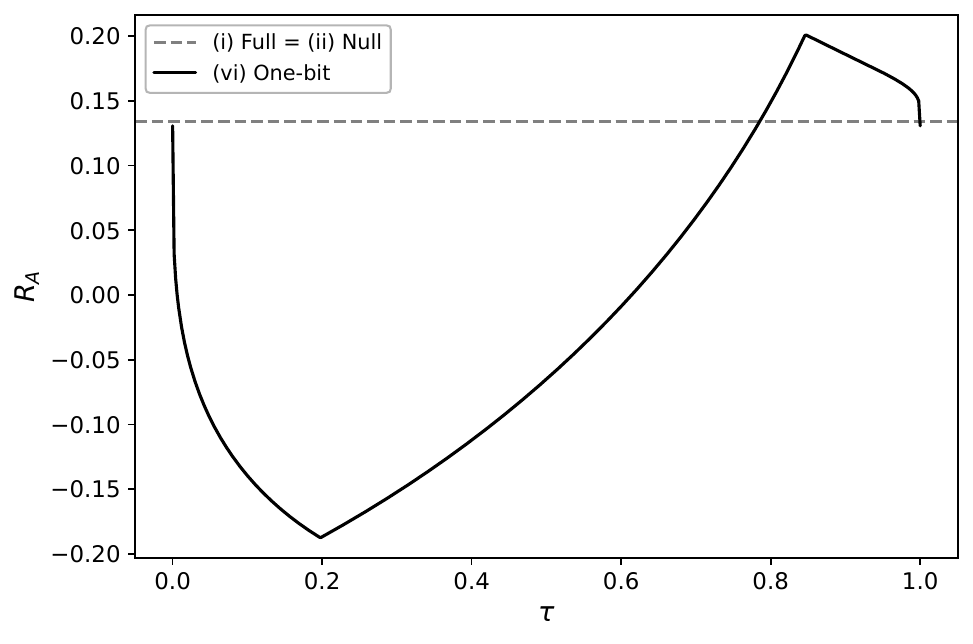}\\{\footnotesize (b) one-bit release vs.\ $\tau$}
    \end{minipage}
    \caption{Integrated risk $R_A$ for the mean adversary ($\sigma_0=1$, $c_E^{\mathrm{mean}}=\sigma_0/2$). \textbf{(a)} noisy full (solid), noisy mean (dash-dotted) and noisy median (dotted) releases as a function of the noise level $\sigma$; \textbf{(b)} the one-bit release $q^{1\mathrm{bit}}_\tau$ as a function of the threshold $\tau$, whose minimum lies off centre at $\tau^\star\approx0.20$ rather than $\tfrac12$. The horizontal grey line is the calibrated reference $R_A(\pi,q^{\dagger})=R_A(\pi,q^{\mathrm{full}})$.}
    \label{fig-gaussian-mean}
\end{figure}

Because of sufficiency and conjugacy, releasing the full vector $X$ with additive noise $\sigma$ or releasing $\bar X$ with noise $\sigma/\sqrt n$ induces identical posterior distributions for $\theta$ and therefore identical integrated risks $R_A$. This equivalence is illustrated in Figure~\ref{fig-gaussian-mean}(a). Minimising $R_A$ over $\sigma$ yields an optimal noise level that balances statistical utility and privacy. (Additional plots displaying the individual risks $R_B$ and $R_E$ for Bob and Eve, respectively, are provided in the appendix.)

Comparing the noisy median release to the noisy full and noisy mean releases, 
we observe that the median achieves the same minimal integrated risk as the mean at the corresponding optimal noise level $\sigma^\star$. 
However, the median release exhibits substantially greater robustness to suboptimal choices of $\sigma$, 
particularly as $\sigma \to 0$, because the exact median is neither sufficient for $\theta$ nor equal to Eve's target $\bar X$.

For the one-bit mechanism, Figure~\ref{fig-gaussian-mean}(b) plots $R_A$ against the single threshold $\tau$; the minimum lies at $\tau^\star\approx0.20\neq\tfrac12$, so releasing Bob's decision (at $\tau=\tfrac12$) is suboptimal here.

After optimising over $\sigma$ and $\tau$ within each class, we compare the resulting risks in Table~\ref{tab-gaussian-mean}. Recall that $\lambda$ is calibrated so that $R_A(\pi,q^{\mathrm{full}})=R_A(\pi,q^{\dagger})$.

\begin{table}[ht]
    \centering
    \begin{tabular}{c|ccc}
        Mechanism & $R_B$ & $R_E$ & $R_A$\\
        \hline
        Full & 0.13 & 0.00 & 0.13 \\
        Null & 0.50 & 0.32 & 0.13 \\
        Noisy full ($\sigma^\star = 3.04$) & 0.31 & 0.25 & 0.03 \\
        Noisy mean ($\sigma^\star = 1.37$) & 0.31 & 0.25 & 0.03 \\
        Noisy median ($\sigma^\star = 1.44$) & 0.32 & 0.26 & 0.03 \\
        One-bit ($\tau^\star = 0.20$) & 0.18 & 0.32 & \textbf{-0.19}
    \end{tabular}
    \caption{Integrated risks $R_B$, $R_E$ and $R_A$ when Eve targets the sample mean at $c_E^{\mathrm{mean}}=\sigma_0/2$, with $\sigma_0=1$ and $\lambda = 1.13$ calibrated so that the null and full releases coincide. Parameters are optimised within each mechanism family. Risks are Monte Carlo estimates rounded to two decimals; $\sigma^\star$ and $\tau^\star$ are grid values.}
    \label{tab-gaussian-mean}
\end{table}

Overall the one-bit release attains the lowest integrated risk, lying strictly below every additive-noise mechanism, as the risk--utility map of Figure~\ref{fig:rumap-ex2}(a) makes plain. Its optimum, however, sits off centre, below $\tfrac12$. At that threshold the released bit is the most informative summary of Bob's problem that still leaves Eve deciding on her prior alone, so Eve's risk stays at its largest possible value and the privacy bound is tight, while Bob's risk has not yet reached its floor, i.e.\ $R_B>R_B^{\mathrm{full}}$. Publishing Bob's decision at $\tau=\tfrac12$ does the opposite, making Bob's utility maximal, i.e.\ $R_B=R_B^{\mathrm{full}}$; but because Eve's mean test is correlated with that decision, revealing said decision lowers her risk below its prior value and loosens the privacy bound. The two bounds cannot be tight at once, so the corner $R_A^{\star}$ of Equation~\eqref{equ-corner-test} is out of reach and the best one-bit release falls short of it.

\subsection{Privacy against the Tails}
\label{sec-gaussian-max}

We now consider the alternative setting where Eve targets the upper tail of the sample, taking 
$T(x)=\max_i x_i$.

In contrast to the previous setting, Bob and Eve are now interested in structurally different features of the data. Bob's decision depends only on $\bar X$, while Eve's decision depends on extreme values. Since $\bar X$ is sufficient for $\theta$ while Eve's target also depends on the ancillary residual $D=\max_i(X_i-\bar X)$, compressing the data towards $\bar X$ conceals the maximum while largely preserving Bob's inference, so a release can help Bob without helping Eve.

\begin{figure}[ht]
    \centering
    \begin{minipage}[b]{0.49\textwidth}\centering
        \includegraphics[width=\textwidth]{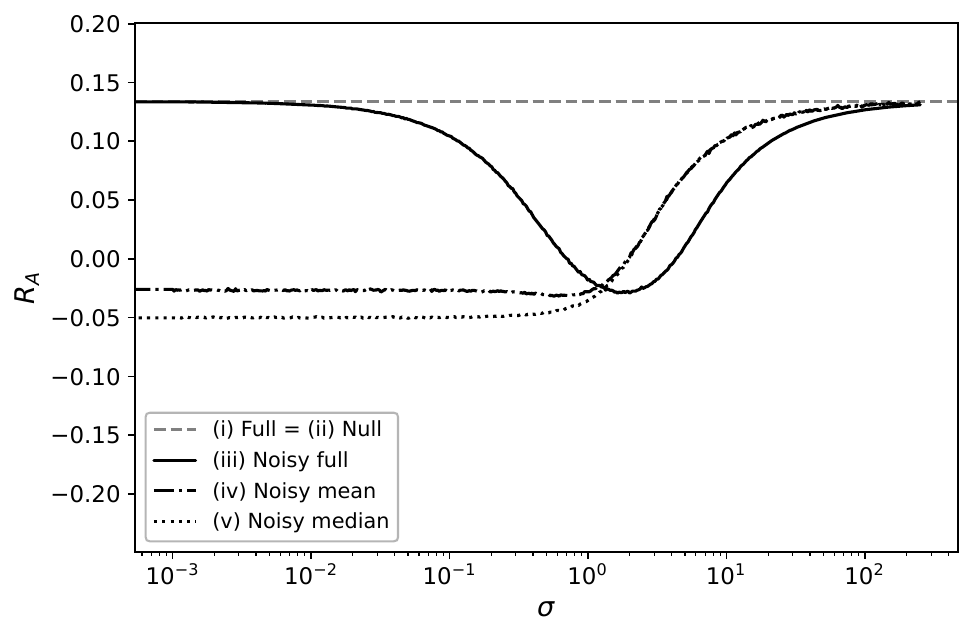}\\{\footnotesize (a) additive-noise releases vs.\ $\sigma$}
    \end{minipage}\hfill
    \begin{minipage}[b]{0.49\textwidth}\centering
        \includegraphics[width=\textwidth]{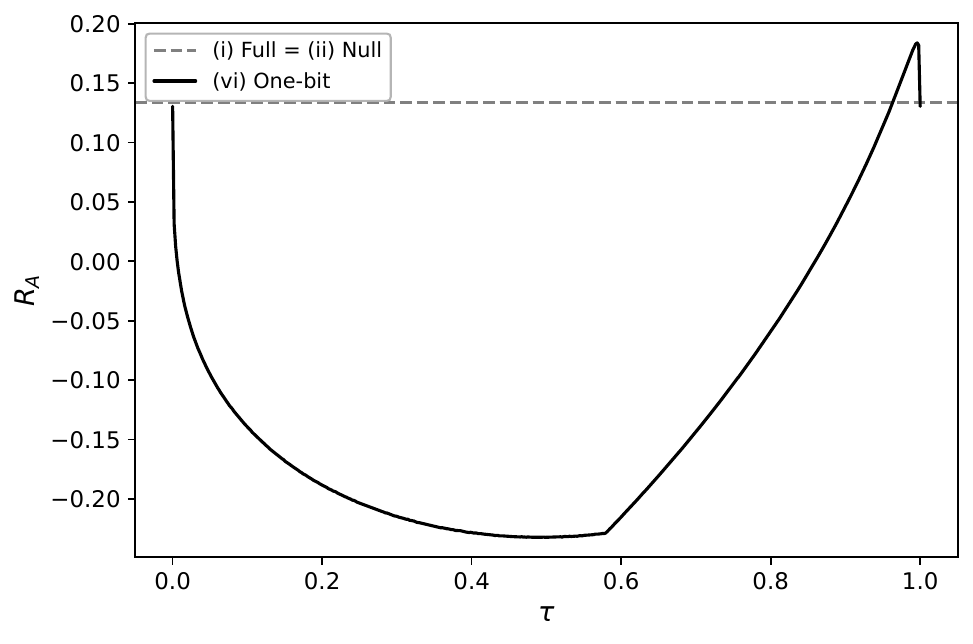}\\{\footnotesize (b) one-bit release vs.\ $\tau$}
    \end{minipage}
    \caption{Integrated risk $R_A$ for the maximum adversary ($\sigma_0=1$, $c_E^{\max}=2\sigma_0$). \textbf{(a)} noisy full (solid), noisy mean (dash-dotted) and noisy median (dotted) releases as a function of the noise level $\sigma$; \textbf{(b)} the one-bit release $q^{1\mathrm{bit}}_\tau$ as a function of the threshold $\tau$, minimised at $\tau^\star=\tfrac12$ where it reaches the corner $R_A\approx-0.23$. The horizontal grey line is the calibrated reference $R_A(\pi,q^{\dagger})=R_A(\pi,q^{\mathrm{full}})$.}
    \label{fig-gaussian-max}
\end{figure}

The additive-noise releases separate along this structural line (Figure~\ref{fig-gaussian-max}(a)). As in the mean case, the noisy full release bridges full disclosure at $\sigma=0$, where $R_E=0$, to the null release as $\sigma\to\infty$, with an interior optimum in between. The noisy mean, however, no longer coincides with full disclosure at $\sigma=0$: it then releases $\bar X$ alone, which this time is not sufficient for Eve's target and so already conceals the maximum, i.e.\ $R_E>0$. Its optimum is therefore reached at a much smaller, and shallow, noise level.

The noisy median release tracks the mean at moderate noise but is strictly increasing in $\sigma$, so its optimum lies at zero noise: the exact median already reveals little about the maximum, and adding noise only costs Bob. It is in fact marginally the best additive release (Table~\ref{tab-gaussian-max}).

In the one-bit release (Figure~\ref{fig-gaussian-max}(b)), the separation between Eve's two goals is starkest, and this time its optimum sits at $\tau=\tfrac12$. Since Alice releases Bob's decision $\eta=\delta_B(x,q^{\mathrm{full}})$ itself, Bob's risk falls to its full-data minimum, $R_B=R_B^{\mathrm{full}}$. Because the two tests depend on different components of the data, and Eve's hypothesis test has significant prior odds, revealing this bit does not change her Bayes decision, only her posterior, so her risk stays at its prior maximum (Table~\ref{tab-gaussian-max}). Both bounds are therefore tight at once, and the one-bit release attains the best possible risk $R_A^{\star}$ of Equation~\eqref{equ-corner-test}, outperforming every other mechanism.

After optimisation within each class, Table~\ref{tab-gaussian-max} compares the resulting risks.

\begin{table}[ht]
    \centering
    \begin{tabular}{c|ccc}
        Mechanism & $R_B$ & $R_E$ & $R_A$\\
        \hline
        Full & 0.13 & 0.00 & 0.13 \\
        Null & 0.50 & 0.24 & 0.13 \\
        Noisy full ($\sigma^\star = 1.63$) & 0.23 & 0.17 & -0.03 \\
        Noisy mean ($\sigma^\star = 0.62$) & 0.21 & 0.16 & -0.03 \\
        Noisy median ($\sigma^\star = 0.00$) & 0.16 & 0.14 & -0.05 \\
        One-bit ($\tau^\star = 0.50$) & 0.13 & 0.24 & \textbf{-0.23}
    \end{tabular}
    \caption{Integrated risks $R_B$, $R_E$ and $R_A$ when Eve targets the sample maximum at the tail threshold $c_E^{\max}=2\sigma_0=2$ ($\sigma_0=1$), with $\lambda=1.52$ calibrated so that the null and full releases coincide. Parameters are optimised within each mechanism family. Risks are Monte Carlo estimates rounded to two decimals; $\sigma^\star$ and $\tau^\star$ are grid values.}
    \label{tab-gaussian-max}
\end{table}

\begin{figure}[ht]
    \centering
    \begin{minipage}{0.49\textwidth}\centering
        \includegraphics[width=\textwidth]{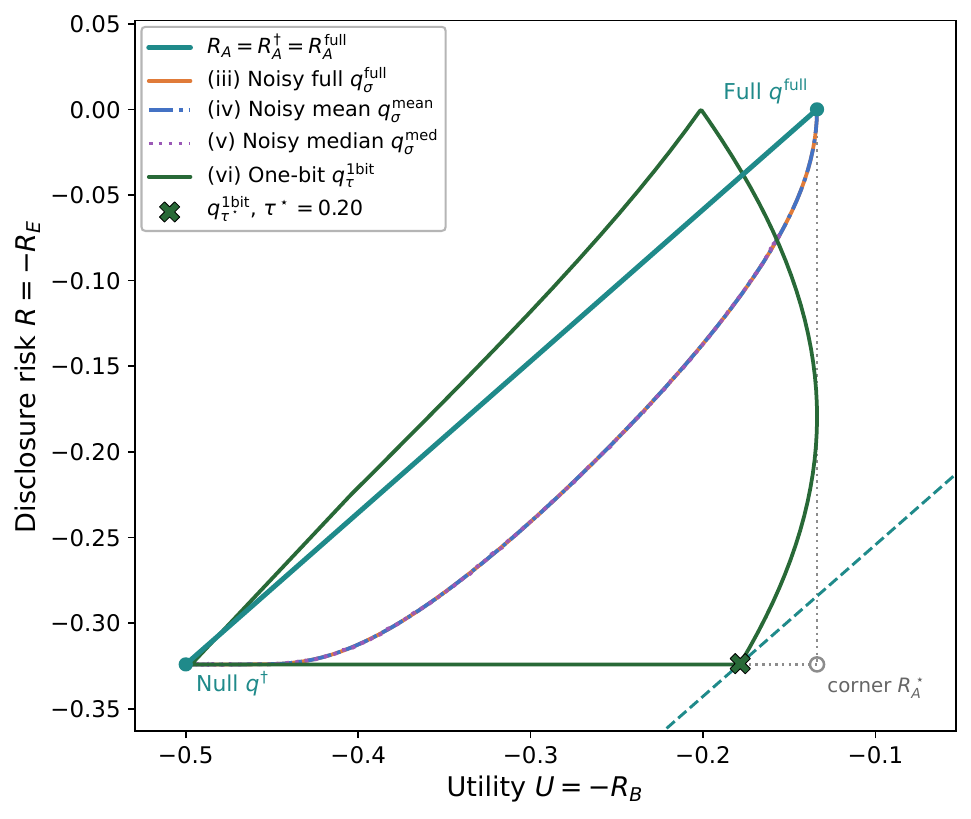}\\{\footnotesize (a) mean adversary ($c_E^{\mathrm{mean}}=\sigma_0/2$)}
    \end{minipage}\hfill
    \begin{minipage}{0.49\textwidth}\centering
        \includegraphics[width=\textwidth]{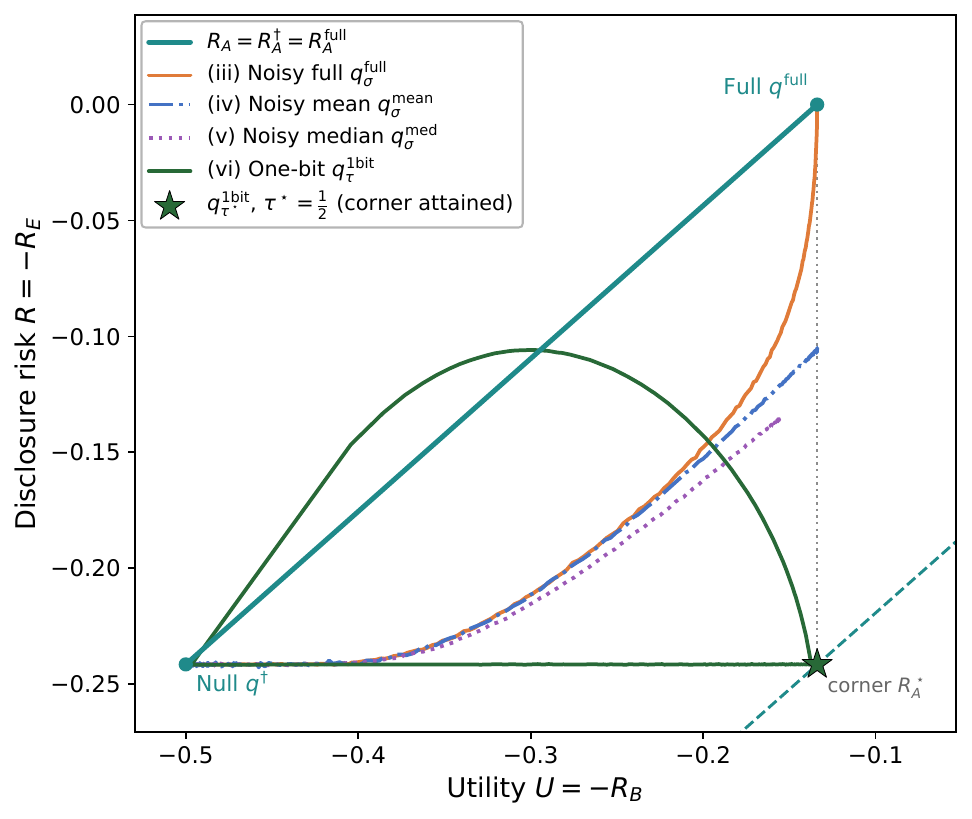}\\{\footnotesize (b) maximum adversary ($c_E^{\max}=2\sigma_0$)}
    \end{minipage}
    \caption{Risk--utility maps for Example~2, read as in Figure~\ref{fig:ru-map}: utility $U=-R_B$ (right, better inference) versus disclosure risk $R=-R_E$ (up, Eve more successful). Full and Null lie on the same level set of $R_A$, and the corner $(-R_B^{\mathrm{full}},-R_E^{\dagger})$, where a vertical from Full meets a horizontal from Null, is the deepest achievable risk $R_A^{\star}$ of Equation~\eqref{equ-corner-test}. \textbf{(a)} Against the mean the corner is \emph{unreachable}: the best release $q^{\star}$ (cross) attains Null's privacy but falls short of Full's utility. \textbf{(b)} Against the maximum the corner is attained by the one-bit release at $\tau=\tfrac12$ (star, $q^{1\mathrm{bit}}_{\tau^\star}=q^{\star}$).}
    \label{fig:rumap-ex2}
\end{figure}

Table~\ref{tab-gaussian-max} exposes this separation: several mechanisms now lower Bob's risk without lowering Eve's, and vice versa, which was impossible when Eve targeted the mean. The risk--utility map of Figure~\ref{fig:rumap-ex2}(b) tells the same story: the one-bit release at $\tau=\tfrac12$ sits exactly at the corner, attaining Full's utility and Null's privacy at once.

\subsection{When Is a Release Optimal?}
\label{sec-gaussian-differences}

Across these scenarios the ranking of the additive-noise releases changes from one case to the next, so no mechanism is uniformly optimal. The one-bit release, a noisy version of Bob's optimal decision $\delta_B(x,q^{\mathrm{full}})$, happens to be best in every case considered here, although it need not be so in general. What is structural rather than incidental is that, among the releases considered, publishing Bob's decision itself, i.e.\ the one-bit release at $\tau=\tfrac12$, is the only one able to reach the optimality corner $R_A^\star$ of Equation~\eqref{equ-corner-test}, where Bob's full-data utility and Eve's null-release privacy hold at once. Whether it does reach the corner is governed by the balance between two forces. What matters is not prior imbalance or task alignment in isolation, but the balance between them: Eve's prior belief determines how firmly she already favours one decision, while the relationship between both tasks determines how much Bob's released decision changes that belief. 

Weak alignment is sufficient to preserve privacy, but it is not necessary: even an informative release may leave Eve's optimal decision unchanged when her prior belief is sufficiently decisive. Conversely, when Eve is a priori close to indifference, even a modest association with Bob's decision may alter her optimal action and improve her inference. This transition is characterised in closed form in the appendix. 
The multiple-Eve criterion of the appendix (Section~S8), where two Eves target the mean and the maximum and Alice is penalised as soon as either is learnt, behaves like the mean adversary rather than the maximum: its mean component stays correlated with Bob's decision, so the corner remains out of reach and the one-bit optimum again sits off centre, at $\tau^\star\approx0.31$. The exact attainment observed against the maximum should nonetheless be read with caution: it rests on the tail test being both a priori decided and only weakly aligned with Bob's decision, and in realistic problems the statistician's and the adversary's targets are rarely so cleanly separated. Furthermore, we can only expect to do this when Eve's decision space is discrete: if she can continuously shift her Bayes rule based on the released information, any non-null release will inevitably change Eve's decision.
The corner is then better understood as a bound on what any release can achieve than as an optimum generally reached.

%% file: sections/conclusion.tex
This paper sets a novel environment for conducting Bayesian statistics under privacy constraints. As such, it allows for the uncertainty quantification that is naturally associated with Bayesian inference, as well as privacy assessment through the posterior distribution of the data conditional on the released quantity. This posterior is quite uncommon and differs from the standard predictive posterior distribution, which is often advocated for model checking \citep{gelman:etal:2013}. However, it complements the original statistical model to induce a complete ranking of release mechanisms. This framework can obviously be expanded in the setting of generalised Bayesian inference \citep{Dempster1968,Bissirietal2016,Matsubaraetal2022}, including score-based methods, variational approximations, or approximate MCMC procedures. 

More centrally,  our privacy framework stands fully embedded within a Bayesian decision-theoretic perspective. While this forces the analyst to add a privacy loss $L_E$ to their Bayesian model, it makes explicit and quantifies which breach of privacy is deemed most crucial. Rather than aiming at a generic, one-size-fits-all privacy protection, like DP, the current approach is tailored for a specific purpose. It is further balanced against the inferential objectives, in an inherent trade-off between disclosure protection and statistical utility. The linear option we adopted offers several benefits, in the spirit of multi-agent decision making, with the coefficient $\lambda$ open to {\em ex post} calibration.

\subsection{Limitations}

Our framework does by no means claim to substitute for existing privacy canons. Instead, it brings a different apprehension of privacy, for which a range of modelling perspectives must be kept in mind. To this end the following may be seen as limitations of our perspective.

Firstly, our framework ranks disclosure decisions averaged over the prior, rather than conditioned on the data. Therefore, the $q$ chosen may fail to optimise Alice's posterior expected loss conditional on the observed data, as explained in Section \ref{section-alice}. However this is inevitable, in order for Eve and Bob to construct well-defined posterior distributions on their respective quantities of interest. This ex ante perspective is also the one found in Bayesian experimental design \citep{Chaloner1995BED}.


Secondly, our framework balances the data harms with data utilities in a utilitarian way. This allows for a single scalar, Alice's risk $R_A(q)$, to summarise the suitability of a release mechanism $q$, and therefore permits the ranking of mechanisms. However it also means that large privacy leaks could in theory be justified for high data utility. More min-max formulations, like that of \citet{bon2026persuasiveprivacy}, may be more comforting for data respondents, but move further away from classical Bayesian analysis parameterised by a loss function. This potential danger of high privacy leaks can be counteracted by choosing losses that are sufficiently sensitive to unwanted data harms.

Thirdly, the manner in which one constructs the Bayesian statistical model for $(x,\theta)$, chooses both loss functions $L_B$ and $L_E$, and calibrates the balancing parameter $\lambda>0$, is by no means obvious. Different choices for these aspects may obviously lead to dissimilar evaluations of release mechanisms through Alice's risk $R_A(q)$. Data utilities and harms can be wide reaching, and many concerns need be weighed up in order to find a single loss function assessing each. However, since it is impossible to simultaneously protect all aspects of the data and perform accurate statistical analysis, the trade-offs between these aspects are fruitful to investigate.

Finally, the quantities introduced in this paper, namely densities, expectations, and losses, may not always be immediately tractable. For this reason computational methods will need be designed towards evaluating Eve and Bob's risk or optimal decisions in complicated enough situations.

\subsection{Extensions}

There are a range of extensions we would like to add to the framework we have introduced in this paper. On a theoretical level, we are interested in investigating the underlying principles that dictate the suitability of a release mechanism $q$. That is to say, we are interested in studying in what situation one may find optimal release mechanisms, and how well these optimal mechanisms perform. For certain classes of losses, for instance quadratic losses for parameter estimation,
a precise construction may be found. For other, more complicated, losses, near-optimal release mechanisms may be found by \textit{discretising} the data space, as well as Bob and Eve's decision spaces, and applying the linear programming strategy introduced in Section \ref{subsec:linear_prog}.

It could be argued that Alice, Eve and Bob's underlying statistical model need not be the same, since they may act under distinct assumptions. However, we would argue that it is sufficient to consider one model for all three, since Bob reflects Alice's interests, and Eve reflects Alice's ``data fears''. In particular if Eve has a different prior distribution this can only perform worse by the standards of the prior distribution Alice selects (if not her decision rule with respect to Alice's prior would not be optimal). Therefore, a more robust perspective may involve instead having families of separate scenarios encompassing different priors, loss functions, values for $\lambda>0$, over which a variety of ex ante assessments $R_1(q),...,R_n(q)$ can be provided. Taking the maximum over this family provides a more conservative assessment of a release mechanism $q$, and similarly a more conservative ranking of mechanisms can be given by a min-max formulation over such a family. This is similar in spirit to the approach of \citet{bon2026persuasiveprivacy}.

From a practical perspective, we are interested in considering a Bayesian formulation of metrics considered in SDC \citep{hundepool2024handbook}. That is to say, it would be useful to try to build prior distributions, statistical models and loss functions that match the data harms and data utility considered in this literature, towards comparing approaches for processing data-sets. This could then allow us to convert the insights gleaned by many decades of successful real-world implementation into a Bayesian framework, providing a posterior predictive model assessment in the spirit of \citet{gelman:etal:2013}. Utilising real-world data would be the final step in this consideration. Similarly, applying these methods in the case of healthcare data used by medical researchers would be of great interest. In both cases, there are two technical challenges we are interested in covering: that of handling relatively high dimensional data rich in co-variates, and that of protecting against \textit{linkage attacks} wherein prior information contains a portion of the full data-set (the difficulties of these two challenges will likely relate to one another).

Another important area of interest arises from the computational requirements of finding and evaluating optimal decisions over Bayesian posteriors conditioned on insufficient statistics \citep{luciano2024insufficient}. Since both Bob and Eve's posteriors require distributions conditionally on $\eta$, these will not generally be available in closed form. It is therefore vital for Alice to be able to accurately calculate and assess their Bayes rules so that she can correctly calculate the risk of a mechanism. This likely will require the development of problem-specific Monte Carlo methods, in the spirit of the data-augmentation MCMC of \citet{ju2022dampc} for privatised summaries, which targets the joint posterior of $(\theta,x)$ and thus delivers Bob's and Eve's posteriors at once.

%% file: supp_body.tex
\section{Setup}

\subsection{Generative model}

We consider the hierarchical Gaussian model
\[
\theta \sim \mathcal{N}(0,\sigma_0^2),
\qquad
X_i \mid \theta \stackrel{\text{i.i.d.}}{\sim} \mathcal{N}(\theta,1),
\quad i=1,\dots,n.
\]

Let $X=(X_1,\dots,X_n)$ denote the dataset and write
\[
\bar X = \frac{1}{n}\sum_{i=1}^n X_i.
\]

A release mechanism is a Markov kernel $q(\dif \eta \mid x)$ on a measurable space $(\mathcal H,\mathcal B)$, producing a released variable $\eta$.

Given $q$, the joint prior predictive distribution of $(\theta,X,\eta)$ factorises as
\begin{equation}
\pi(\dif \theta)\, p(\dif x\mid\theta)\, q(\dif \eta\mid x).
\label{eq:joint_theta_x_eta}
\end{equation}

All integrated risks are evaluated with respect to this joint distribution.

Throughout, $\Phi$ denotes the standard normal cumulative distribution function and $\bar\Phi=1-\Phi$ its survival function.

\subsection{Decision problems and Bayes $0$--$1$ risk}
\label{sec:decision}

Bob and Eve each face a binary decision problem under $0$--$1$ loss.

\paragraph{Bob (utility objective).}
Bob tests the parameter-level hypotheses
\[
H_{0,B} : \theta \le c_B,
\qquad
H_{1,B} : \theta > c_B,
\]
and we define the latent state indicator
\[
S_B = \ind\{\theta > c_B\}.
\]

\paragraph{Eve (privacy objective).}
Eve tests data-level hypotheses
\[
H_{0,E} : T(X) \le c_E,
\qquad
H_{1,E} : T(X) > c_E,
\]
where $T(X)$ is a data-functional. In this appendix, we consider
\[
T(X)=\bar X
\quad \text{(mean adversary)},
\qquad
T(X)=\max_{1\le i\le n} X_i
\quad \text{(maximum adversary)}.
\]
The corresponding state indicator is
\[
S_E = \ind\{T(X) > c_E\}.
\]

\paragraph{Integrated risks.} Both $S_B$ and $S_E$ are unobserved at the decision stage and must be inferred from the released variable $\eta$. For each agent $a \in \{B,E\}$, define the posterior success probability
\begin{equation}
p_a(\eta) = \Prob(S_a = 1 \mid \eta).
\label{eq:posterior_success_prob}
\end{equation}

Under $0$--$1$ loss, each agent adopts the Bayes optimal decision rule, 
i.e.\ chooses the most probable state given $\eta$. 
The resulting minimal conditional Bayes risk is therefore
\[
r_a(\eta)
=
r\!\bigl(p_a(\eta)\bigr)
=
\min\{p_a(\eta),\,1-p_a(\eta)\},
\]
where $r(p)=\min\{p,1-p\}$ denotes the optimal $0$--$1$ Bayes risk function.

The integrated risks are
\begin{equation}
    R_a(\pi,q) = \E[r_a(\eta)],
    \label{eq:integrated_risk}
\end{equation}

where the expectation is taken with respect to the joint distribution defined in Equation~\eqref{eq:joint_theta_x_eta}.

In this work, we define Alice’s \emph{ex ante} risk as a linear combination of the two agents’ integrated Bayes risks. The parameter $\lambda>0$ governs the privacy–utility trade-off by penalising Eve’s performance relative to Bob’s:

\begin{equation*}
R_A(\pi,q)
=
R_B(\pi,q)
-
\lambda R_E(\pi,q).
\end{equation*}

The risk is termed \emph{ex ante} because it is evaluated under the joint prior predictive distribution, prior to observing any specific dataset.

In the following sections, we show how the agents can compute their respective risks of interest under the different release mechanisms. Closed-form expressions are rarely available, so we describe for each case the numerical procedure used to evaluate these quantities and to produce the tables and figures reported in the paper.

\section{Case 1: Full release}
\label{sec:full}
\subsection{Mechanism}

As a first benchmark, we consider the \emph{full release} mechanism, denoted $q^{\mathrm{full}}$, under which Alice releases the entire data vector. That is,
\[
\eta = X.
\]

\subsection{Bob: computation of $R_B(\pi,q^{\mathrm{full}})$}

Under full release, Bob observes $X$ and therefore works with the standard posterior distribution, which is conjugate in the Gaussian model:
\begin{equation}
\theta \mid X \sim \mathcal{N}(m_X,v_X),\quad\text{with}\quad v_X = (\sigma_0^{-2}+n)^{-1}
\text{ and }
m_X = v_X n \bar X.
\label{eq:standard_posterior}
\end{equation}

The posterior success probability is then known exactly as
\[
p_B(X)
=
\Prob(\theta>c_B \mid X)
=
\bar\Phi\!\left(
\frac{c_B-m_X}{\sqrt{v_X}}
\right).
\]


Importantly, $p_B(X)$ depends on the data only through the sample mean $\bar X$. 
Since, under the prior predictive distribution,
\begin{equation}
\bar X \sim \mathcal{N}\!\left(0,v_{\bar X}\right), \qquad v_{\bar X} = \sigma_0^2 + \frac 1 n
\label{eq:prior_pred_barx}
\end{equation}
the integrated risk reduces to a one-dimensional expectation:
\[
R_B(\pi,q^{\mathrm{full}})
=
\E_{\bar X}
\!\left[
\min\!\left(
\bar\Phi\!\left(
\frac{c_B - v_X n \bar X}{\sqrt{v_X}}
\right),
\Phi\!\left(
\frac{c_B - v_X n \bar X}{\sqrt{v_X}}
\right)
\right)
\right].
\]

This integral is one-dimensional and Gaussian, and can therefore be evaluated efficiently using standard quadrature methods, such as Gauss--Hermite quadrature.

\subsection{Eve: computation of $R_E(\pi,q^{\mathrm{full}})$}

Under full release, Eve observes the entire dataset $X$. 
In both cases $T(X)=\bar X$ (mean adversary) and 
$T(X)=\max_{1\le i\le n}X_i$ (maximum adversary), the target statistic is observed exactly. 
The event $\{T(X)>c_E\}$ is therefore known without uncertainty, implying
\[
r_E(X)=0.
\]
Consequently,
\[
R_E(\pi,q^{\mathrm{full}})=0.
\]


\section{Case 2: Null release}

\subsection{Mechanism}

Under the \emph{null release} mechanism, denoted $q^\dagger$, no information about the data is disclosed. The mechanism outputs a constant, independently of the data $X$,
\[
\eta = \dagger.
\]

Consequently, observing $\eta$ provides no information about either $\theta$ or $X$, and the posterior distribution coincides with the prior.

\subsection{Bob: computation of $R_B(\pi,q^\dagger)$}

Since no information is released, Bob must rely entirely on his prior beliefs, namely
\[
\theta \sim \mathcal{N}(0,\sigma_0^2).
\]

His posterior success probability therefore coincides with the prior probability:
\[
p_B(\dagger)
=
\Prob(\theta>c_B)
=
\bar\Phi\!\left(\frac{c_B}{\sigma_0}\right).
\]

As this quantity is deterministic and independent of the data, the integrated risk is available in closed form:
\[
R_B(\pi,q^\dagger)
=
\min\!\left\{
\bar\Phi\!\left(\frac{c_B}{\sigma_0}\right),
\Phi\!\left(\frac{c_B}{\sigma_0}\right)
\right\}.
\]

\subsection{Eve: computation of $R_E(\pi,q^\dagger)$}

Similarly, Eve must rely entirely on the prior predictive distribution of the target statistic $T(X)$.

\paragraph{Mean adversary.}

When $T(X)=\bar X$, Eve can compute the success probability under the prior predictive distribution defined in Equation~\eqref{eq:prior_pred_barx}: 
\[
p_E(\dagger)
=
\Prob(\bar X>c_E)
=
\bar\Phi\!\left(
\frac{c_E}{\sqrt{v_{\bar X}}}
\right),
\]
which again is deterministic. The integrated risk is therefore
\[
R_E(\pi,q^\dagger)
=
\min\!\left\{
\bar\Phi\!\left(
\frac{c_E}{\sqrt{v_{\bar X}}}
\right),
\Phi\!\left(
\frac{c_E}{\sqrt{v_{\bar X}}}
\right)
\right\}.
\]

\paragraph{Maximum adversary.}

When $T(X)=\max_{1\le i\le n} X_i$, the $X_i$ are independent conditionally on $\theta$, with
\[
X_i \mid \theta \sim \mathcal{N}(\theta,1).
\]

Thus, conditionally on $\theta$,
\[
\Prob\!\left(
\max_{1\le i\le n} X_i \le c_E
\mid \theta
\right)
=
\Phi(c_E-\theta)^n.
\]

Integrating over the prior $\theta \sim \mathcal{N}(0,\sigma_0^2)$ yields
\[
p_E(\dagger)
=
1-
\E_\theta
\!\left[
\Phi(c_E-\theta)^n
\right].
\]

This is a one-dimensional Gaussian expectation and can be evaluated numerically (e.g.\ via Gauss--Hermite quadrature). The integrated risk is then
\[
R_E(\pi,q^\dagger)
=
r\!\left(p_E(\dagger)\right).
\]

\section{Case 3: Noisy full release}
\subsection{Mechanism}

After the two extreme mechanisms (full and null release), a natural intermediate mechanism is the \emph{noisy full release}. Under this mechanism, independent Gaussian noise with standard deviation $\sigma \ge 0$ is added to each observation. We denote it $q_\sigma^{\mathrm{full}}$, and the released variable is $\eta=Y$, where
\[
Y_i = X_i + \varepsilon_i,
\qquad
\varepsilon_i \stackrel{\mathrm{i.i.d.}}{\sim} \mathcal{N}(0,\sigma^2),
\]
independently of $(\theta,X)$. Conditionally on $\theta$, we then have
\[
Y_i \mid \theta \sim \mathcal{N}(\theta,1+\sigma^2).
\]
\subsection{Bob: computation of $R_B(\pi,q_\sigma^{\mathrm{full}})$}

Since the augmented model remains Gaussian, conjugacy still holds. Given the release $\eta=Y$, Bob can compute the posterior distribution in closed form:
\begin{equation}
\begin{gathered}
\theta \mid Y \sim \mathcal{N}(m_\sigma^{\mathrm{full}}(\bar Y), v_\sigma^{\mathrm{full}}),\\[4pt]
\text{with}\quad v_\sigma^{\mathrm{full}}
=
\left(
\sigma_0^{-2}
+
\frac{n}{1+\sigma^2}
\right)^{-1}, \quad 
m_\sigma^{\mathrm{full}}(\bar Y)
=
v_\sigma^{\mathrm{full}}
\frac{n}{1+\sigma^2}
\,\bar Y.
\end{gathered}
\label{eq:noisy_posterior}
\end{equation}

The posterior success probability is therefore
\[
p_B(Y)
=
\bar\Phi
\!\left(
\frac{c_B-m_\sigma^{\mathrm{full}}(\bar Y)}{\sqrt{v_\sigma^{\mathrm{full}}}}
\right).
\]

As in the full release case, $p_B(Y)$ depends on the data only through $\bar Y$.
Under the prior predictive distribution,
\begin{equation}
\bar Y
\sim
\mathcal{N}
\!\left(
0,
v_{\bar Y}
\right),\qquad v_{\bar Y} = \sigma_0^2 + \frac{1+\sigma^2}{n}.
\label{eq:prior_pred_bary}
\end{equation}

Bob’s integrated risk therefore reduces to a one-dimensional Gaussian expectation:
\[
R_B(\pi,q_\sigma^{\mathrm{full}})
=
\E_{\bar Y}
\!\left[
\min\!\left(
\bar\Phi\!\left(
\frac{c_B-m_\sigma^{\mathrm{full}}(\bar Y)}{\sqrt{v_\sigma^{\mathrm{full}}}}
\right),
\Phi\!\left(
\frac{c_B-m_\sigma^{\mathrm{full}}(\bar Y)}{\sqrt{v_\sigma^{\mathrm{full}}}}
\right)
\right)
\right],
\]
which can be evaluated numerically using standard quadrature methods.

\subsection{Eve: computation of $R_E(\pi,q_\sigma^{\mathrm{full}})$}

\paragraph{Mean adversary.}

Eve aims to infer the event $\{\bar X>c_E\}$. Under the noisy release,
\[
\bar Y=\bar X+\bar\varepsilon,
\qquad
\bar\varepsilon\sim\mathcal{N}\!\left(0,\frac{\sigma^2}{n}\right),
\]
with $\bar{\varepsilon}$ independent of $\bar X$.

Since $(\bar X,\bar Y)$ are jointly Gaussian, conditioning yields
\[
\bar X \mid \bar Y
\sim
\mathcal{N}\!\left(a\,\bar Y,\;v_{\bar X}(1-a)\right),
\qquad
a=\frac{v_{\bar X}}{v_{\bar X}+\sigma^2/n},
\]
with $v_{\bar X}$ defined in Equation~\eqref{eq:prior_pred_barx}. Hence, the probability of interest is known in closed form
\[
p_E(\bar Y)
=
\bar\Phi\!\left(
\frac{c_E-a\bar Y}
{\sqrt{v_{\bar X}(1-a)}}
\right).
\]

Therefore, the integrated risk can be evaluated numerically via quadrature, since it reduces to an expectation with respect to the one-dimensional Gaussian prior predictive distribution defined in Equation~\eqref{eq:prior_pred_bary}, i.e.
\[
R_E(\pi,q_\sigma^{\mathrm{full}})
=
\E_{\bar Y}
\!\left[
\min\{p_E(\bar Y),1-p_E(\bar Y)\}
\right].
\]

\paragraph{Maximum adversary.}

Eve now seeks to infer the event $\{\max_{1\le i\le n} X_i>c_E\}$. 
Applying Bayes’ rule together with the posterior distribution of $\theta \mid Y$ given in Equation~\eqref{eq:noisy_posterior}, we obtain
\[
p_E(Y)
=
1-
\Prob\!\left(
\max_i X_i \le c_E
\mid Y
\right)
=
1-
\E_{\theta \mid Y}
\!\left[
\Prob\!\left(
\max_i X_i \le c_E
\mid \theta,Y
\right)
\right].
\]

Conditionally on $(\theta,Y)$, the coordinates $X_i$ are independent Gaussian variables. More precisely,
\[
X_i \mid (\theta,Y_i)
\sim
\mathcal{N}\!\left(\mu_i(\theta,Y_i),s^2\right),
\]
where
\[
\mu_i(\theta,Y_i)
=
Y_i+\frac{\sigma^2}{1+\sigma^2}(\theta-Y_i),
\qquad
s^2=\frac{\sigma^2}{1+\sigma^2}.
\]

It follows that, conditionally on $(\theta,Y)$,
\[
\Prob\!\left(
\max_i X_i \le c_E
\mid \theta,Y
\right)
=
\prod_{i=1}^n
\Phi\!\left(
\frac{c_E-\mu_i(\theta,Y_i)}{s}
\right).
\]

Substituting this expression into the previous identity yields
\[
p_E(Y)
=
1-
\E_{\theta \mid Y}
\!\left[
\prod_{i=1}^n
\Phi\!\left(
\frac{c_E-\mu_i(\theta,Y_i)}{s}
\right)
\right].
\]

Finally, the integrated risk is obtained by averaging over the prior predictive distribution of $Y$:
\[
R_E(\pi,q_\sigma^{\mathrm{full}})
=
\E_Y\!\left[r_E(Y)\right].
\]

The inner expectation with respect to $\theta\mid Y$ is one-dimensional and is evaluated via Gauss--Hermite quadrature, while the outer expectation over $Y\in\mathbb{R}^n$ is computed by Monte Carlo under the prior predictive $Y \sim \mathcal{N}\!\bigl(\mathbf{0}_n,\ \sigma_0^2\,\mathbf{1}_n\mathbf{1}_n^{\!\top} + (1+\sigma^2)\, I_n\bigr)$, whose coordinates are correlated through the shared $\theta$. In practice the code samples $\theta$, then $X\mid\theta$, then $Y\mid X$ hierarchically, so this correlation is respected.

\section{Case 4: Noisy mean release}

\subsection{Mechanism}

In the same spirit as the previous section, we consider the \emph{noisy mean release}, denoted $q_\sigma^{\mathrm{mean}}$, under which
\[
\eta = \bar X + \varepsilon,
\qquad
\varepsilon \sim \mathcal{N}(0,\sigma^2),
\]
with the noise independent of $(\theta,X)$.

Conditionally on $\theta$, we have
\[
\eta \mid \theta
\sim
\mathcal{N}(\theta,v_\eta),
\qquad
v_\eta = \frac{1}{n}+\sigma^2.
\]

Under the prior predictive distribution,
\begin{equation}
\eta \sim \mathcal{N}(0,\sigma_0^2+v_\eta).
\label{eq:prior_pred_noisy_mean}
\end{equation}

\subsection{Bob: computation of $R_B(\pi,q_\sigma^{\mathrm{mean}})$}

Since the model remains Gaussian, conjugacy holds as in Equation~\eqref{eq:standard_posterior}. Conditioning on the released statistic $\eta$ yields
\begin{equation*}
\theta \mid \eta
\sim
\mathcal{N}(m_\sigma^{\mathrm{mean}}(\eta),v_\sigma^{\mathrm{mean}}),
\label{eq:noisy_mean_posterior}
\end{equation*}
where
\[
v_\sigma^{\mathrm{mean}}
=
\left(
\sigma_0^{-2}+v_\eta^{-1}
\right)^{-1},
\qquad
m_\sigma^{\mathrm{mean}}(\eta)
=
v_\sigma^{\mathrm{mean}}\,
\frac{\eta}{v_\eta}.
\]

Bob tests the event $\{\theta>c_B\}$. His posterior success probability is therefore
\[
p_B(\eta)
=
\bar\Phi\!\left(
\frac{c_B-m_\sigma^{\mathrm{mean}}(\eta)}
{\sqrt{v_\sigma^{\mathrm{mean}}}}
\right).
\]

Since $\eta$ is Gaussian under the prior predictive distribution (Equation~\eqref{eq:prior_pred_noisy_mean}), the integrated risk reduces to the one-dimensional expectation
\[
R_B(\pi,q_\sigma^{\mathrm{mean}})
=
\E_\eta
\!\left[
r_B(\eta)
\right],
\]
which is evaluated numerically via Gauss--Hermite quadrature.

\subsection{Eve: computation of $R_E(\pi,q_\sigma^{\mathrm{mean}})$}

Eve considers two possible targets: the sample mean and the sample maximum.

\paragraph{Mean adversary.}

Eve seeks to infer the event $\{\bar X>c_E\}$. 
Under the noisy mean release, $(\bar X,\eta)$ are jointly Gaussian. 
Using the prior predictive distribution of $\bar X$ (Equation~\eqref{eq:prior_pred_barx}) together with standard Gaussian conditioning, we obtain
\[
\bar X \mid \eta
\sim
\mathcal{N}\!\left(\tilde a\,\eta,\; v_{\bar X}(1-\tilde a)\right),
\]
where
\[
\tilde a
=
\frac{v_{\bar X}}{v_{\bar X}+\sigma^2},
\]
and $v_{\bar X}$ is defined in Equation~\eqref{eq:prior_pred_barx}.

It follows that the posterior success probability is
\[
p_E(\eta)
=
\bar\Phi\!\left(
\frac{c_E-\tilde a\,\eta}
{\sqrt{v_{\bar X}(1-\tilde a)}}
\right).
\]

Since $\eta$ is Gaussian under the prior predictive distribution (Equation~\eqref{eq:prior_pred_noisy_mean}), the integrated risk reduces to the one-dimensional expectation
\[
R_E(\pi,q_\sigma^{\mathrm{mean}})
=
\E_\eta
\!\left[
\min\!\{p_E(\eta),1-p_E(\eta)\}
\right],
\]
which is evaluated numerically via Gauss--Hermite quadrature.

\paragraph{Maximum adversary.}

Eve now considers the event $\{\max_{1\le i\le n} X_i>c_E\}$. A natural but \emph{incorrect}
approach is to marginalise over $\theta$, writing
$\Prob(\max_i X_i\le c_E\mid\eta)=\E_{\theta\mid\eta}[\Phi(c_E-\theta)^n]$: this treats the $X_i$ as
free $\mathcal N(\theta,1)$ variables given $\eta$, whereas the noisy mean release constrains their
average $\bar X$, so that $\max_i X_i$ remains correlated with $\eta$ even after conditioning on
$\theta$. (The identity $\Prob(\max_i X_i\le c_E\mid\theta)=\Phi(c_E-\theta)^n$ holds, but propagating it
through $\theta\mid\eta$ over-estimates $R_E$; for $\sigma\to0$ it gives $0.127$ instead of the correct
$0.079$ at $n=5,\sigma_0=1,c_E=0$, a balanced check distinct from the main-text tail adversary's $c_E=2\sigma_0$.)

The release $\eta=\bar X+\varepsilon$ depends on the data only through $\bar X$, so we condition on
$\bar X$ rather than $\theta$. Decompose
\[
\max_{1\le i\le n} X_i=\bar X+D,\qquad D=\max_{1\le i\le n}\bigl(X_i-\bar X\bigr).
\]
Conditionally on $\theta$, $X_i-\bar X=e_i-\bar e$ with $e_i\stackrel{\text{iid}}{\sim}\mathcal N(0,1)$;
since the vector of deviations $(e_i-\bar e)_i$ is independent of $\bar e$, the statistic $D$ is
independent of $(\bar X,\theta,\eta)$, and its law depends neither on $\theta$ nor on $\sigma$.
Using the mean-adversary conditional law $\bar X\mid\eta\sim\mathcal N\!\bigl(\tilde a\,\eta,\;
v_{\bar X}(1-\tilde a)\bigr)$ derived above, we obtain the exact posterior success probability
\[
p_E(\eta)
=
\Prob\!\left(\bar X+D>c_E\mid\eta\right)
=
1-\E_{D}\!\left[
\Phi\!\left(\frac{c_E-D-\tilde a\,\eta}{\sqrt{v_{\bar X}(1-\tilde a)}}\right)
\right],
\]
and the integrated risk
\[
R_E(\pi,q_\sigma^{\mathrm{mean}})
=
\E_\eta
\!\left[
\min\!\{p_E(\eta),1-p_E(\eta)\}
\right].
\]
In practice the ancillary statistic $D$ is simulated once by Monte Carlo ($M$ draws of
$W\sim\mathcal N(\mathbf 0_n,I_n)$, $D=\max_i(W_i-\bar W)$); the inner expectation over $D$ is
estimated by an empirical mean, while the outer expectation over $\eta$ (one-dimensional Gaussian) is
evaluated by Gauss--Hermite quadrature.

\section{Case 5: Noisy median release}

\subsection{Mechanism}

We consider the \emph{noisy median release} mechanism $q_\sigma^{\mathrm{med}}$, defined by
\[
\eta
=
\mathrm{med}(X)
+
\varepsilon,
\qquad
\varepsilon\sim\mathcal{N}(0,\sigma^2),
\]
with the noise independent of $(\theta,X)$.

Under this mechanism, the conditional distributions
\[
\theta\mid\eta,
\qquad
\bar X\mid\eta,
\qquad
\max_i X_i\mid\eta
\]
do not admit closed-form expressions. We therefore approximate the integrated risks by Monte Carlo.

\subsection{Monte Carlo approximation of the risks}


As defined in Section~\ref{sec:decision}, the posterior success probability $p_a(\eta)$ and the integrated Bayes risk $R_a(\pi,q)$, for $a \in \{B,E\}$, are given in Equations~\eqref{eq:posterior_success_prob} and \eqref{eq:integrated_risk}. Under the noisy median mechanism, closed-form expressions are not available. We therefore resort to Monte Carlo approximation under the joint prior predictive distribution.

\paragraph{Prior predictive simulation.}

We generate Monte Carlo samples
\[
(\theta^{(m)},X^{(m)},\eta^{(m)})_{m=1}^M
\]
as follows:
\[
\theta^{(m)}\sim\pi(\cdot),
\qquad
X^{(m)}\sim p(\cdot\mid\theta^{(m)}),
\qquad
\eta^{(m)}=\mathrm{med}(X^{(m)})+\varepsilon^{(m)},
\]
with $\varepsilon^{(m)}\sim\mathcal{N}(0,\sigma^2)$ independent.

For each agent, we construct the corresponding indicators
\[
S_B^{(m)}=\ind\{\theta^{(m)}>c_B\},
\]
and
\[
S_E^{(m)}=\ind\{T(X^{(m)})>c_E\},
\]
where $T(X)$ is either $\bar X$ (mean adversary) or $\max_i X_i$ (maximum adversary).

\paragraph{Nonparametric estimation of $p_a(\eta)$.} We approximate the conditional success probability
\[
p_a(\eta)=\Prob(S_a=1\mid\eta)
\]
by a histogram-based regression estimator. We construct a uniform grid of $K$ bins on
\[
[-L\sqrt{V_\eta},\, L\sqrt{V_\eta}],
\qquad
V_\eta = \sigma_0^2 + v_{\mathrm{med}} + \sigma^2,
\qquad
v_{\mathrm{med}}=\frac{\pi}{2n},
\]
where $v_{\mathrm{med}}$ is the standard normal approximation to $\Var(\mathrm{med}(X)\mid\theta)$.
Let $(I_k)_{k=1}^K$ denote the resulting bins. For $\eta\in\mathbb R$, let $I_k(\eta)$ be the bin containing $\eta$ and define
\[
\widehat p_a(\eta)
=
\frac{
\sum_{m:\,\eta^{(m)}\in I_k(\eta)} S_a^{(m)} + \alpha
}{
\sum_{m:\,\eta^{(m)}\in I_k(\eta)} 1 + 2\alpha
},
\qquad
\alpha=\tfrac12,
\]
where $\alpha$ corresponds to Jeffreys smoothing (a $\mathrm{Beta}(1/2,1/2)$ prior) to stabilise the estimate in low-count bins.
In all experiments we use $L=6$, $K=200$ and $M=4\times10^{6}$ prior predictive draws. 

\paragraph{Approximation of the integrated risk.}

The integrated risk is then approximated by
\[
\widehat R_a
=
\frac{1}{M}
\sum_{m=1}^M
r\!\left(
\widehat p_a(\eta^{(m)})
\right),
\]
which provides a Monte Carlo approximation of
\[
R_a(\pi,q_\sigma^{\mathrm{med}})
=
\E[r_a(\eta)].
\]

This procedure is applied identically for Bob and for Eve, the only difference lying in the definition of the event $S_a$.

\section{Case 6: One-bit release}

\subsection{Mechanism}

We now consider a \emph{one-bit release} mechanism, denoted $q_\tau^{1\mathrm{bit}}$, defined by
\[
\eta
=
\ind\{p_B^{\mathrm{full}}(\bar X) > \tau\},
\]
where $p_B^{\mathrm{full}}(\bar X)$ is Bob’s posterior success probability in the full release case (see Section~\ref{sec:full}), computed from the standard posterior in Equation~\eqref{eq:standard_posterior}:
\[
p_B^{\mathrm{full}}(\bar X)
=
\Prob(\theta>c_B \mid \bar X)
=
\bar\Phi\!\left(
\frac{c_B - v_X n \bar X}{\sqrt{v_X}}
\right).
\]

Since the mapping $\bar X \mapsto p_B^{\mathrm{full}}(\bar X)$ is monotone, the mechanism is equivalent to thresholding $\bar X$:
\[
\eta
=
\ind\{\bar X > t_\tau\},
\]
where the deterministic threshold is
$
t_\tau=
\frac{
c_B - \sqrt{v_X}\,\Phi^{-1}(1-\tau)
}{
n\,v_X
}$ with $v_X = (\sigma_0^{-2} + n)^{-1}$ as defined in Equation~\eqref{eq:standard_posterior}.

\subsection{Bob: computation of $R_B(\pi,q_\tau^{1\mathrm{bit}})$}

Since $\eta\in\{0,1\}$, the integrated risk decomposes as
\[
R_B(\pi,q_\tau^{1\mathrm{bit}})
=
\sum_{j=0}^1
\Prob(\eta=j)\,
r\!\left(
\Prob(\theta>c_B\mid\eta=j)
\right).
\]

All required quantities are expressed as expectations with respect to the prior on $\theta$.

Conditionally on $\theta$, we have
\[
\Prob(\eta=1\mid\theta)
=
\Prob(\bar X>t_\tau\mid\theta)
=
\bar\Phi\!\left(\sqrt n\,(t_\tau-\theta)\right).
\]

Therefore,
\[
\Prob(\eta=1)
=
\E_\theta
\!\left[
\bar\Phi\!\left(\sqrt n\,(t_\tau-\theta)\right)
\right],
\qquad
\theta\sim\mathcal N(0,\sigma_0^2).
\]

For $j=1$,
\[
\Prob(\theta>c_B,\eta=1)
=
\E_\theta
\!\left[
\ind\{\theta>c_B\}
\bar\Phi\!\left(\sqrt n\,(t_\tau-\theta)\right)
\right],
\]
and for $j=0$,
\[
\Prob(\theta>c_B,\eta=0)
=
\E_\theta
\!\left[
\ind\{\theta>c_B\}
\Phi\!\left(\sqrt n\,(t_\tau-\theta)\right)
\right].
\]

The posterior probabilities are then obtained by Bayes’ rule:
\[
\Prob(\theta>c_B\mid\eta=j)
=
\frac{
\Prob(\theta>c_B,\eta=j)
}{
\Prob(\eta=j)
}.
\]

All expectations above are one-dimensional Gaussian integrals in $\theta$. 
Following the implementation, we evaluate them via Gauss--Hermite quadrature after the change of variable $\theta=\sigma_0 z$.

\subsection{Eve: computation of $R_E(\pi,q_\tau^{1\mathrm{bit}})$}

\paragraph{Mean adversary.}

Since $\eta=\ind\{\bar X>t_\tau\}$ and $\bar X\sim\mathcal N(0,v_{\bar X})$ under the prior predictive distribution (Equation~\eqref{eq:prior_pred_barx}), the conditional law $\bar X\mid\eta$ is truncated normal. Hence,
\[
\Prob(\bar X>c_E\mid \eta=1)
=
\Prob(\bar X>c_E\mid \bar X>t_\tau)
=
\frac{\Prob\!\left(\bar X>\max\{c_E,t_\tau\}\right)}{\Prob(\bar X>t_\tau)}
=
\frac{\bar\Phi\!\left(\frac{\max\{c_E,t_\tau\}}{\sqrt{v_{\bar X}}}\right)}
{\bar\Phi\!\left(\frac{t_\tau}{\sqrt{v_{\bar X}}}\right)},
\]
and
\[
\Prob(\bar X>c_E\mid \eta=0)
=
\Prob(\bar X>c_E\mid \bar X\le t_\tau)
=
\frac{\Prob(c_E<\bar X\le t_\tau)}{\Prob(\bar X\le t_\tau)}
=
\frac{\bigl[\Phi\!\left(\frac{t_\tau}{\sqrt{v_{\bar X}}}\right)-\Phi\!\left(\frac{c_E}{\sqrt{v_{\bar X}}}\right)\bigr]_+}
{\Phi\!\left(\frac{t_\tau}{\sqrt{v_{\bar X}}}\right)},
\]
where $[x]_+=\max\{x,0\}$.

The integrated risk decomposes as
\[
R_E(\pi,q_\tau^{1\mathrm{bit}})
=
\sum_{j=0}^1
\Prob(\eta=j)\,
r\!\left(
\Prob(\bar X>c_E \mid \eta=j)
\right),
\]
and all terms are available in closed form.

\paragraph{Maximum adversary.}

Conditioning on $\eta$ restricts $\bar X$ but does not yield a tractable expression for the joint distribution of $X$. 
We therefore approximate
\[
\Prob(\max_i X_i>c_E \mid \eta=j)
\]
by Monte Carlo under the joint prior predictive distribution.

Specifically, we simulate $M$ independent draws $(\theta^{(m)},X^{(m)})$ from the prior predictive model, compute $\bar X^{(m)}$, define $\eta^{(m)}=\ind\{\bar X^{(m)}>t_\tau\}$, and evaluate the indicator
\[
S_E^{(m)}=\ind\{\max_i X_i^{(m)}>c_E\}.
\]
The conditional probabilities are then estimated as empirical averages within the subsets $\{\eta=j\}$:
\[
\widehat{\Prob}(\max_i X_i>c_E \mid \eta=j)
=
\frac{
\sum_{m:\,\eta^{(m)}=j} S_E^{(m)}
}{
\sum_{m:\,\eta^{(m)}=j} 1
}.
\]

The integrated risk is then approximated by
\[
\widehat R_E
=
\sum_{j=0}^1
\widehat{\Prob}(\eta=j)\,
r\!\left(
\widehat{\Prob}(\max_i X_i>c_E \mid \eta=j)
\right).
\]

\section{Case 7: Multiple Eves}
\label{sec:supp-multiple-eves}

\subsection{Mechanism and criterion}

So far Eve has pursued a single feature. Suppose instead that Alice wishes to protect \emph{several} features simultaneously, so that a breach occurs as soon as Eve learns \emph{any one} of them. Following the ``multiple Bobs or Eves'' remark of the main text (Remark~3.1), we combine two Eves, one testing the sample mean ($T_1=\ind\{\bar X>c_1\}$) and one testing the maximum ($T_2=\ind\{\max_i X_i>c_2\}$), through the \emph{minimum} of their losses, $L_E=\min(L_E^{1},L_E^{2})$, the alternative aggregation mentioned in that remark, so that Eve incurs no loss as long as she is correct about at least one target. We call this the \emph{multiple-Eve criterion}. Since Eve is penalised only when wrong about \emph{both} targets, her Bayes-optimal play is the complement of the least-probable joint state (not the joint MAP assignment), and the integrated privacy risk is exactly that smallest joint-cell probability,
\begin{equation}
R_E(\pi,q)
=
\E_\eta\!\left[\min_{t_1,t_2\in\{0,1\}} P_{t_1t_2}(\eta)\right],
\qquad
P_{t_1t_2}(\eta)=\Prob\!\left(T_1{=}t_1,\,T_2{=}t_2\mid\eta\right),
\label{eq:supp-neither}
\end{equation}
that is, $R_E=\Prob(\text{Eve is wrong about both})$: Alice succeeds only if Eve learns neither feature. Bob's problem is unchanged, and the same six release mechanisms as in Cases~1--6 are compared.

\subsection{Choice of thresholds}

Criterion~\eqref{eq:supp-neither} is degenerate when the two events are a priori \emph{nested}: since $\bar X\le\max_i X_i$, a common threshold $c_1=c_2$ gives $\{\bar X>c\}\subseteq\{\max_i X_i>c\}$, the state $(1,0)$ is impossible, Eve plays the straddle $(\delta_1,\delta_2)=(0,1)$ and always secures one answer, so that $R_E\equiv0$. It also collapses when the outlier threshold is too extreme, the state $(0,1)=\{\bar X\le c_1,\max_i X_i>c_2\}$ then requiring an atypically large deviation and becoming negligible. The criterion is informative only for a \emph{moderate} outlier threshold at which both off-diagonal states keep positive mass; we reuse the single-adversary thresholds $c_1=\sigma_0/2$ and $c_2=2\sigma_0$ with $\sigma_0=1$, the choice $c_1>0$ preserving the state $(0,1)$ (a larger $\sigma_0$ would re-nest the events). Even then $R_E$ remains small, reflecting that ``Eve learns neither'' is an intrinsically demanding requirement; a more robust way to control several adversaries at once is to take a worst-case minimum over the individual risks, or a weighted sum of losses, as in Remark~3.1 of the main text.

\subsection{Computation of $R_E$}

\paragraph{Joint cells via the mean/deviation decomposition.}
Reusing $\max_i X_i=\bar X+D$ with $D$ ancillary (independent of $\bar X,\theta,\eta$; see the maximum
adversary of Case~4), the two targets become thresholds on the \emph{same} variable
$\bar X\mid\eta\sim\mathcal N(u,s^2)$, at levels $c_1$ and $c_2-D$. Writing
$G(t)=\Phi\!\left((t-u)/s\right)$, $\ell=\min(c_1,c_2-D)$, $h=\max(c_1,c_2-D)$,
\begin{align*}
P_{00}&=\E_D[G(\ell)], & P_{11}&=\E_D[1-G(h)],\\
P_{10}&=\E_D\bigl[(G(h)-G(\ell))\,\ind{c_2-D>c_1}\bigr], &
P_{01}&=\E_D\bigl[(G(h)-G(\ell))\,\ind{c_2-D\le c_1}\bigr].
\end{align*}
For the null and noisy mean releases $(u,s^2)$ is the (prior or posterior) Gaussian law of $\bar X$
and the cells are semi-analytic; for the one-bit release $\bar X\mid\eta$ is a truncated normal; for
the noisy median and noisy full releases, where $D$ is no longer ancillary given $\eta$, the cells
are estimated by Monte Carlo (binning in $\eta$, respectively an inner conditional simulation of
$X\mid Y$).

\subsection{Results}

\begin{table}[ht]
    \centering
    \begin{tabular}{c|ccc}
        Mechanism & $R_B$ & $R_E$ & $R_A$\\
        \hline
        Full & 0.13 & 0.00 & 0.13 \\
        Null & 0.50 & 0.03 & 0.13 \\
        Noisy full ($\sigma^\star = 1.80$) & 0.24 & 0.03 & -0.10 \\
        Noisy mean ($\sigma^\star = 0.63$) & 0.21 & 0.03 & -0.11 \\
        Noisy median ($\sigma^\star = 0.85$) & 0.25 & 0.03 & -0.08 \\
        One-bit ($\tau^\star = 0.31$) & 0.15 & 0.03 & \textbf{-0.18}
    \end{tabular}
    \caption[Integrated risks under the multiple-Eve criterion]{Integrated risks under the multiple-Eve criterion of Case~7, with $\sigma_0=1$, $c_1=\sigma_0/2$, $c_2=2\sigma_0$ and $\lambda\approx12.4$. Since $R_E$ is small for every release, $R_A$ is dominated by the term $\lambda R_E$ and its ranking, rather than its values, is the meaningful output; see the text.}
    \label{tab-gaussian-neither}
\end{table}

\begin{figure}[ht]
    \centering
    \begin{minipage}[b]{0.49\textwidth}\centering
        \includegraphics[width=\textwidth]{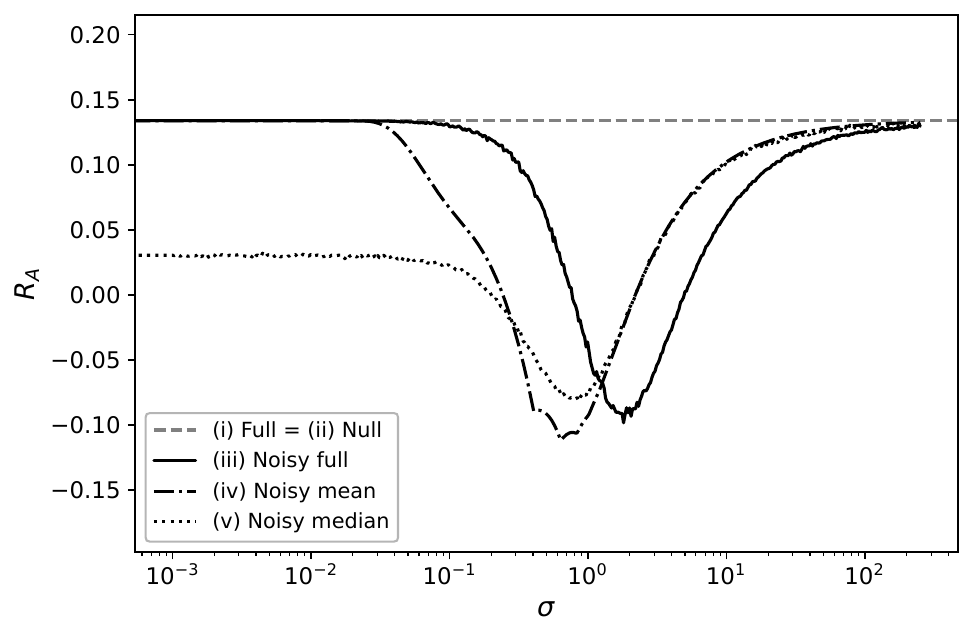}\\{\footnotesize (a) additive-noise releases vs.\ $\sigma$}
    \end{minipage}\hfill
    \begin{minipage}[b]{0.49\textwidth}\centering
        \includegraphics[width=\textwidth]{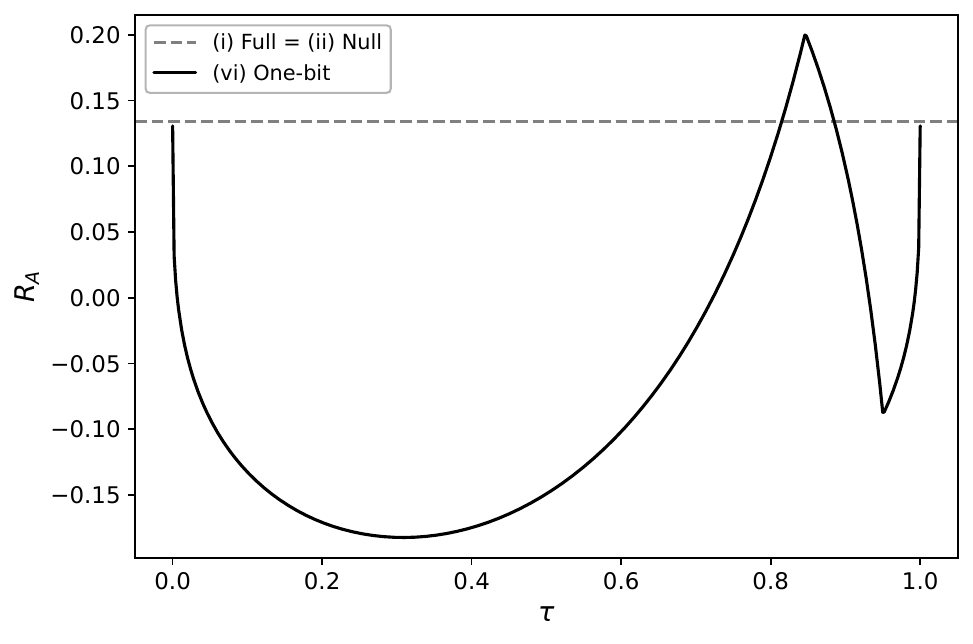}\\{\footnotesize (b) one-bit release vs.\ $\tau$}
    \end{minipage}
    \caption[Integrated risk for the multiple-Eve criterion]{Integrated risk $R_A$ for the multiple-Eve criterion ($\sigma_0=1$). \textbf{(a)} noisy full/mean/median releases as a function of $\sigma$; \textbf{(b)} the one-bit release as a function of $\tau$. The grey line is the calibrated reference $R_A(\pi,q^{\dagger})=R_A(\pi,q^{\mathrm{full}})$.}
    \label{fig-gaussian-neither}
\end{figure}

Two points stand out. First, $R_E$ stays \emph{small} ($\lesssim0.03$) even under the null release: the two events nearly straddle, so Eve can almost always guess at least one feature and ``learning neither'' is intrinsically hard to enforce, hence the large calibrated $\lambda$ and the qualitative reading of $R_A$. Second, the one-bit release again attains the lowest risk, as in both single-feature cases, though the additive releases now rank differently; its optimum is again off-centre, i.e.\ $\tau^\star\approx0.31$: as for the pure mean adversary, the mean component of the joint criterion is aligned with Bob's decision, so releasing it at $\tau=\tfrac12$ stays informative to Eve. 
Section~\ref{sec:supp-onebit-opt} explains why the corner of the risk--utility map is lost here, as for the mean adversary, and Section~\ref{sec:supp-rumap} reproduces the corresponding map.

\section{Summary of integrated risk computations}

For clarity, we summarise here the numerical strategies used to evaluate the integrated risks for each release mechanism and adversary.
In all cases, expectations are taken with respect to the joint prior predictive distribution induced by $(\pi,p,q)$.

Depending on the structure of the mechanism, the computation reduces to either closed-form expressions, one-dimensional Gaussian integrals evaluated via Gauss--Hermite quadrature (GH), nested quadrature, or Monte Carlo (MC) simulation under the prior predictive model.

Table~\ref{tab:methods} summarises, for each release mechanism and single-feature adversary, 
the key conditional probability driving the Bayes risk together with the numerical method used for evaluation. The multiple-Eve criterion of Case~7 reuses the same building blocks, the mean-adversary conditional law of $\bar X$ and the ancillary deviation $D$, as detailed in Section~\ref{sec:supp-multiple-eves}.

\begin{table}[h]
\centering
\begin{tabular}{lll}
\hline
Mechanism & Agent $a$ & Evaluation of $R_a(\pi,q)$ \\
\hline

Full release 
& Bob 
& $\E_{\bar X}\!\left[r\!\left(\Prob(\theta>c_B\mid \bar X)\right)\right]$ \\
& Eve (mean/max) 
& $0$ \\[4pt]

Null release 
& Bob 
& $r\!\left(\Prob(\theta>c_B)\right)$ (exact) \\
& Eve (mean) 
& $r\!\left(\Prob(\bar X>c_E)\right)$ (exact) \\
& Eve (max) 
& $r\!\left(1-\E_\theta[\Phi(c_E-\theta)^n]\right)$ (1D GH) \\[4pt]

Noisy full release 
& Bob 
& $\E_{\bar Y}\!\left[r\!\left(\Prob(\theta>c_B\mid \bar Y)\right)\right]$ (1D GH) \\
& Eve (mean) 
& $\E_{\bar Y}\!\left[r\!\left(\Prob(\bar X>c_E\mid \bar Y)\right)\right]$ (1D GH) \\
& Eve (max) 
& $\E_Y\!\left[r\!\left(\Prob(\max_i X_i>c_E\mid Y)\right)\right]$ (MC + 1D GH) \\[4pt]

Noisy mean release 
& Bob 
& $\E_{\eta}\!\left[r\!\left(\Prob(\theta>c_B\mid \eta)\right)\right]$ (1D GH) \\
& Eve (mean) 
& $\E_{\eta}\!\left[r\!\left(\Prob(\bar X>c_E\mid \eta)\right)\right]$ (1D GH) \\
& Eve (max)
& $\E_{\eta}\!\left[r\!\left(1-\E_{D}\bigl[\Phi\bigl(\tfrac{c_E-D-\tilde a\eta}{\sqrt{v_{\bar X}(1-\tilde a)}}\bigr)\bigr]\right)\right]$ (MC) \\[4pt]

Noisy median release
& Bob/Eve 
& $\E\!\left[r\!\left(\Prob(S_a=1\mid \eta)\right)\right]$ (MC + histogram regression) \\[4pt]

One-bit release 
& Bob 
& $\sum_{j=0}^1 \Prob(\eta=j)\,
r\!\left(\Prob(\theta>c_B\mid \eta=j)\right)$ (1D GH) \\
& Eve (mean) 
& discrete truncated-normal formulas (exact) \\
& Eve (max) 
& discrete decomposition + MC \\

\hline
\end{tabular}
\caption{Summary of the integrated risk computations for each release mechanism and adversary. 
All expectations are taken under the prior predictive distribution induced by $(\pi,p,q)$. 
GH denotes Gauss--Hermite quadrature; MC denotes Monte Carlo simulation.}
\label{tab:methods}
\end{table}

\section{When does the one-bit release attain the global optimum?}
\label{sec:supp-onebit-opt}

Across every adversary considered above the one-bit release performs well, and at the balanced
threshold $\tau=\tfrac12$ it often minimises $R_A$ over \emph{all} mechanisms, not merely within the
one-bit family. We give here the underlying reason, together with a closed-form characterisation of when it holds. Throughout, $\pi=\Prob(\theta>c_B)=\bar\Phi(c_B/\sigma_0)$ is Bob's prior
success probability and $p=\Prob(T(X)>c_E)$ is Eve's.

\paragraph{Two structural bounds.}
Because every release $\eta$ is a (possibly randomised) function of the data, Bob's posterior on
$\theta$ is a garbling of the full-data posterior; as $0$--$1$ Bayes risk cannot increase under
garbling,
\[
R_B(\pi,q)\ \ge\ R_B(\pi,q^{\mathrm{full}})\ =:\ R_B^{\mathrm{full}}.
\]
Conversely, writing $\kappa(\eta)=\Prob(T(X)>c_E\mid\eta)$ for Eve's posterior, the tower rule
gives $\E[\kappa(\eta)]=p$, and since $r(u)=\min(u,1-u)$ is concave, Jensen's inequality yields
\[
R_E(\pi,q)=\E\bigl[r(\kappa(\eta))\bigr]\ \le\ r(p)=\min(p,1-p).
\]
Bob's risk is bounded \emph{below} by the full-release value; Eve's is bounded \emph{above} by her
prior risk.

\paragraph{A universal corner.}
Recall the calibration $\lambda=(R_B^{\dagger}-R_B^{\mathrm{full}})/\min(p,1-p)$, where
$R_B^{\dagger}=\min(\pi,1-\pi)$ is Bob's null-release risk (Case~2), chosen so the null and full
releases share the same $R_A$. Substituting the two bounds into $R_A=R_B-\lambda R_E$,
\[
R_A(\pi,q)\ \ge\ R_A^\star\ :=\ 2R_B^{\mathrm{full}}-\min(\pi,1-\pi),
\]
a lower bound met with equality only when \emph{both} bounds are tight simultaneously, a corner
of the risk--utility map, whose depth is governed by Bob's full-release risk
$R_B^{\mathrm{full}}$ together with the balance of his prior test $\min(\pi,1-\pi)$.

\paragraph{When the one-bit release reaches the corner.}
At $\tau=\tfrac12$ the released bit is exactly Bob's Bayes decision under full disclosure,
$\eta=\delta_B=\ind{p_B^{\mathrm{full}}(\bar X)>\tfrac12}$. This decision is a sufficient reduction for
Bob's $0$--$1$ problem, so the first bound holds with equality, $R_B=R_B^{\mathrm{full}}$ exactly. The
second bound is tight iff Jensen is tight over the two-point law of $\kappa$, that is iff the two
posteriors $\kappa_0,\kappa_1$ (the values of $\kappa(\eta)$ at $\eta=0,1$) lie on the \emph{same} side
of $\tfrac12$: seeing $\delta_B$ never flips Eve's decision, so $R_E$ stays at its prior value
$\min(p,1-p)$. This is the \emph{plateau} regime, and there the one-bit release at $\tau=\tfrac12$ is
globally optimal, attaining $R_A^\star$.

\paragraph{Closed-form condition.}
Let $\beta=\Prob(\delta_B=1)$ be Bob's decision rate and $s=\kappa_1-\kappa_0>0$ (positive
association between the two tests). Since $p=\beta\kappa_1+(1-\beta)\kappa_0$, the posteriors straddle
$\tfrac12$, and the corner is \emph{lost}, exactly on a band of $p$,
\[
\text{corner reached}\iff p\notin\Bigl(\tfrac12-(1-\beta)s,\ \tfrac12+\beta s\Bigr).
\]
The band has width $s$ and centre $\tfrac12+(\beta-\tfrac12)s$; it is centred at $\tfrac12$ only
when Bob's test is balanced ($\beta=\tfrac12$).

\paragraph{The mean adversary in closed form.}
For $T=\bar X$ every quantity is explicit: the decision rate is $\beta=\bar\Phi\!\bigl(c_B\sqrt{v_{\bar X}}/\sigma_0^2\bigr)$, and the corner condition reduces to
\[
\text{corner reached}\iff p\le\tfrac{\beta}{2}\quad\text{or}\quad p\ge\tfrac{1+\beta}{2},
\]
a lost-corner band $\bigl(\tfrac{\beta}{2},\tfrac{1+\beta}{2}\bigr)$ of width exactly $\tfrac12$.
Equivalently, holding $p$ fixed, the plateau occupies a half-line in $\tau$ whose edge is exact for
any threshold $c_B$,
\[
\tau_{\min}(p)=\Phi\!\Bigl(\sqrt n\,\sigma_0\,\Phi^{-1}\!\bigl(2(1-p)\bigr)
-\tfrac{c_B\sqrt{n\sigma_0^2+1}}{\sigma_0}\Bigr)\quad(p>\tfrac12),
\]
plateau for $\tau\ge\tau_{\min}$, with the symmetric branch $\tau_{\max}(p)=
\Phi\bigl(\sqrt n\,\sigma_0\,\Phi^{-1}(1-2p)-c_B\sqrt{n\sigma_0^2+1}/\sigma_0\bigr)$ for $p<\tfrac12$.
A non-balanced Bob test enters only through the additive shift
$-\,c_B\sqrt{n\sigma_0^2+1}/\sigma_0$ (recovering $\Phi(\sqrt n\,\sigma_0\,\Phi^{-1}(2(1-p)))$ at
$c_B=0$); imposing $\tau_{\min}(p)=\tfrac12$ returns exactly the band edge $p=\tfrac{1+\beta}{2}$ above,
so the same $c_B$ sets both the position of the band (through $\beta$) and the shift of the
$\tau$-plateau.

\paragraph{Summary.}
The one-bit release publishing Bob's decision is globally optimal precisely when Eve's target is a
priori \emph{lopsided} ($p$ far from $\tfrac12$), the failure band, where revealing $\delta_B$ does
move Eve, being centred near $\tfrac12$ with width $s$ (exactly $\tfrac12$ for the mean). Two forces
act in the same direction: the more \emph{balanced} Bob's test ($\pi\approx\tfrac12$, deepening the
attainable corner $R_A^\star=2R_B^{\mathrm{full}}-\min(\pi,1-\pi)$) and the more a priori
\emph{decided} Eve's test ($p$ away from $\tfrac12$), the more the simple threshold rule at
$\tau^\star=\tfrac12$ is not merely good but exactly optimal. In the main text, the maximum adversary
at $c_E^{\max}=2\sigma_0$ has $p\approx0.24$; its association $s$, and hence its exact band, is available only numerically, but $p$ lies outside it, so releasing $\delta_B$ at
$\tau=\tfrac12$ is globally optimal; the mean adversary at $c_E^{\mathrm{mean}}=\sigma_0/2$ has
$p\approx0.32$, \emph{inside} the band, so its one-bit optimum drifts off $\tau=\tfrac12$ (to
$\tau^\star\approx0.20$). The multiple-Eve criterion inherits this same regime: its response to the released bit is governed by the mean coordinate, whose $p$ lies inside the band, so the corner is likewise lost and its one-bit optimum drifts off $\tau=\tfrac12$ (to $\tau^\star\approx0.31$).

\section{Additional numerical illustrations}

This section provides additional graphical diagnostics for the integrated risks introduced in the main text. 
Whereas the main paper reports the integrated risk $R_A$ only, we additionally display the individual risks of Bob ($R_B$) and Eve ($R_E$) to illustrate the structural mechanisms underlying the privacy–inference trade-offs.

As in the main text, all expectations are evaluated under the joint prior predictive distribution.

\subsection{Variation with the noise level $\sigma$}

\begin{figure}[h]
\centering
\includegraphics[width=\textwidth]{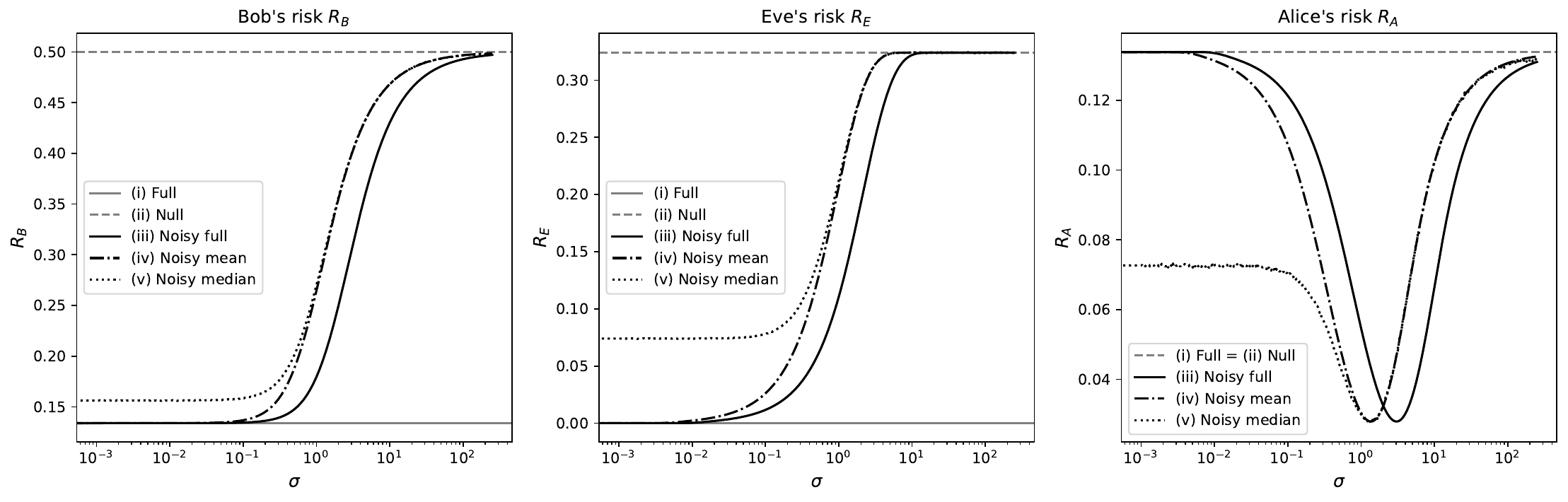}
\caption{Integrated risks as a function of $\sigma$ when Eve targets the sample mean.}
\label{fig:supp-sigma-mean}
\end{figure}

\begin{figure}[h]
\centering
\includegraphics[width=\textwidth]{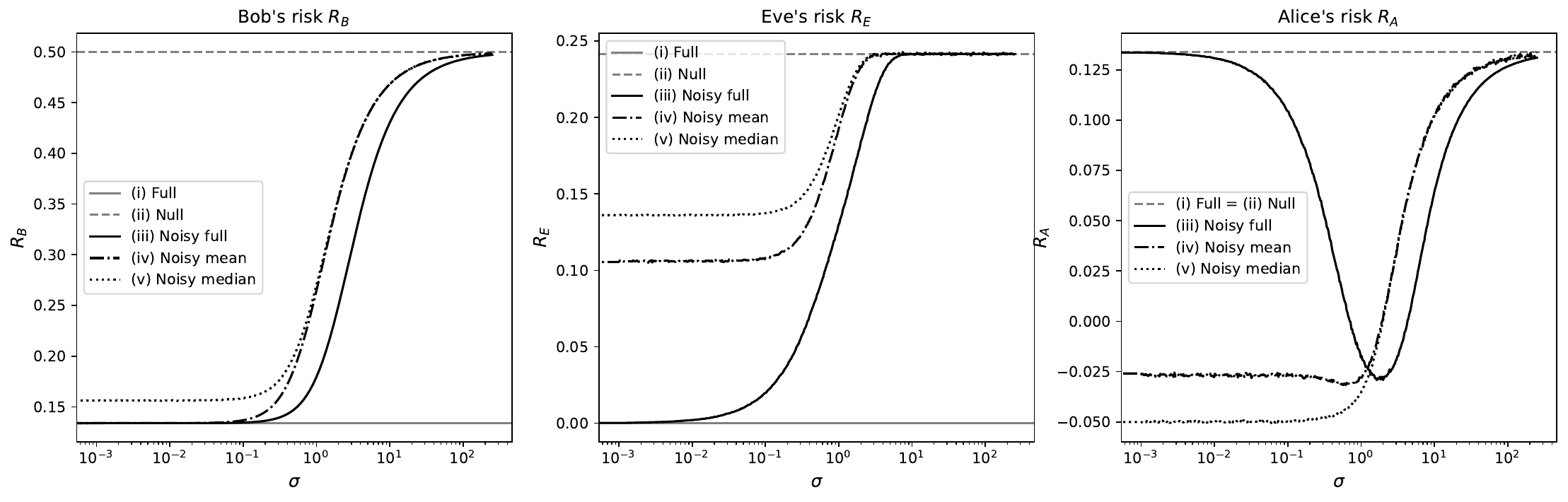}
\caption{Integrated risks as a function of $\sigma$ when Eve targets the sample maximum.}
\label{fig:supp-sigma-max}
\end{figure}

\begin{figure}[h]
\centering
\includegraphics[width=\textwidth]{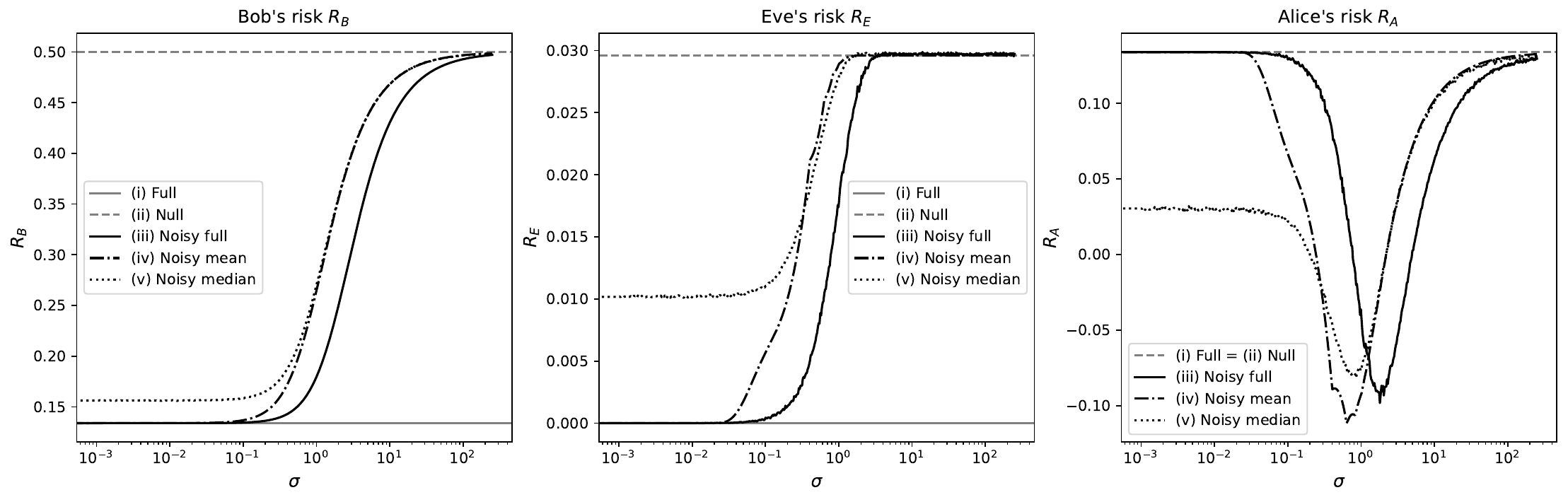}
\caption{Integrated risks $R_B$, $R_E$ and $R_A$ as a function of $\sigma$ under the multiple-Eve criterion.}
\label{fig:supp-sigma-neither}
\end{figure}

Figures~\ref{fig:supp-sigma-mean} and~\ref{fig:supp-sigma-max} show the risks for Cases~3--5 as functions of the noise parameter $\sigma$, for the mean and maximum adversaries respectively, with the multiple-Eve criterion in Figure~\ref{fig:supp-sigma-neither}. 

Each figure consists of three panels corresponding to $R_B$, $R_E$, and $R_A$. 
Cases~1 (full release) and~2 (null release) are shown as horizontal benchmarks.

\subsection{Variation with the threshold $\tau$}

\begin{figure}[h]
\centering
\includegraphics[width=\textwidth]{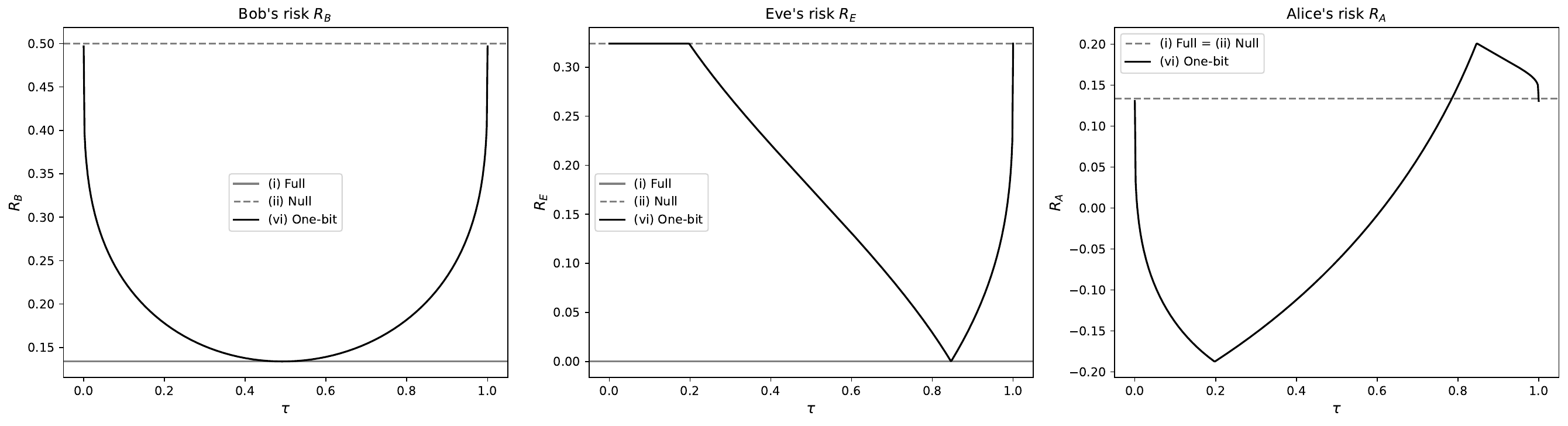}
\caption{Integrated risks as a function of $\tau$ when Eve targets the sample mean.}
\label{fig:supp-tau-mean}
\end{figure}

\begin{figure}[h]
\centering
\includegraphics[width=\textwidth]{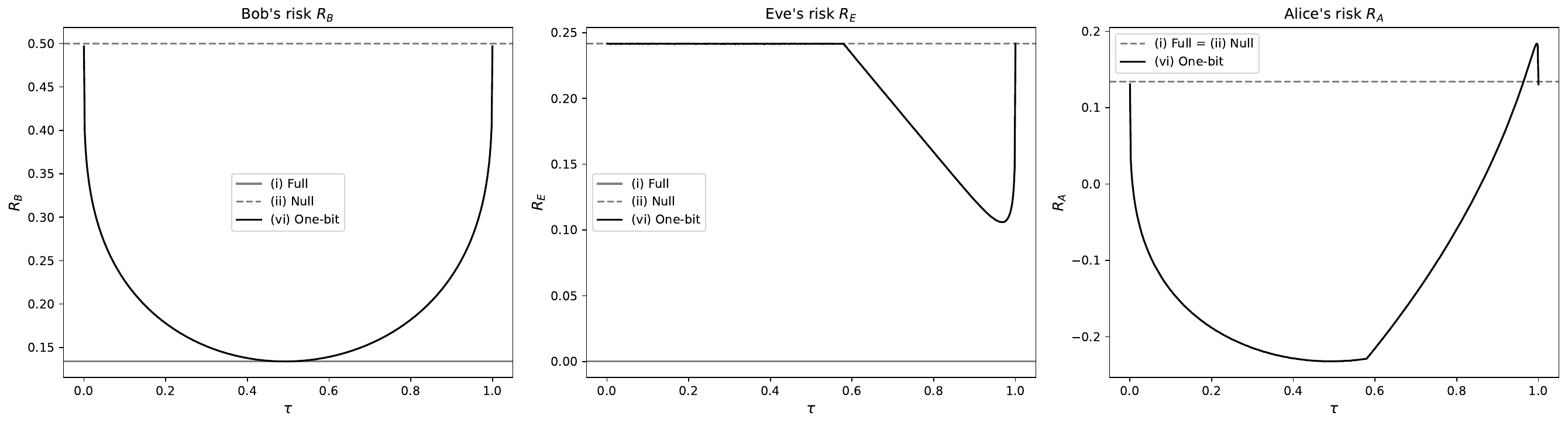}
\caption{Integrated risks as a function of $\tau$ when Eve targets the sample maximum.}
\label{fig:supp-tau-max}
\end{figure}

\begin{figure}[h]
\centering
\includegraphics[width=\textwidth]{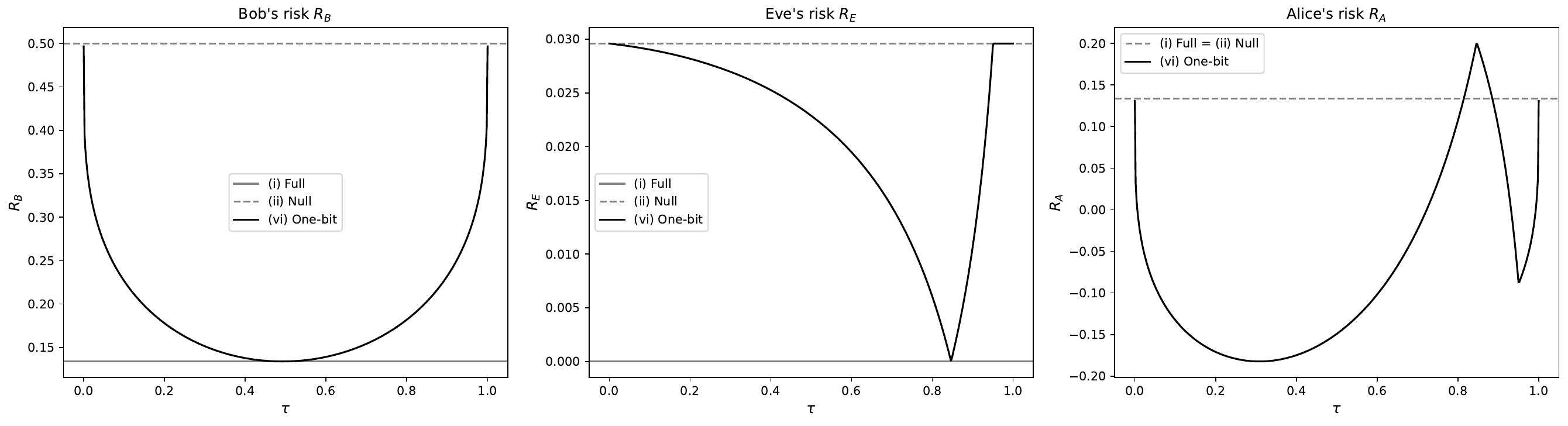}
\caption{Integrated risks $R_B$, $R_E$ and $R_A$ as a function of $\tau$ under the multiple-Eve criterion.}
\label{fig:supp-tau-neither}
\end{figure}

Figures~\ref{fig:supp-tau-mean} and~\ref{fig:supp-tau-max} display the risks for the one-bit mechanism (Case~6) as functions of the threshold parameter $\tau$, again for both adversaries, with the multiple-Eve criterion in Figure~\ref{fig:supp-tau-neither}.

\subsection{Risk--utility maps}
\label{sec:supp-rumap}

\begin{figure}[h]
\centering
\includegraphics[width=0.7\textwidth]{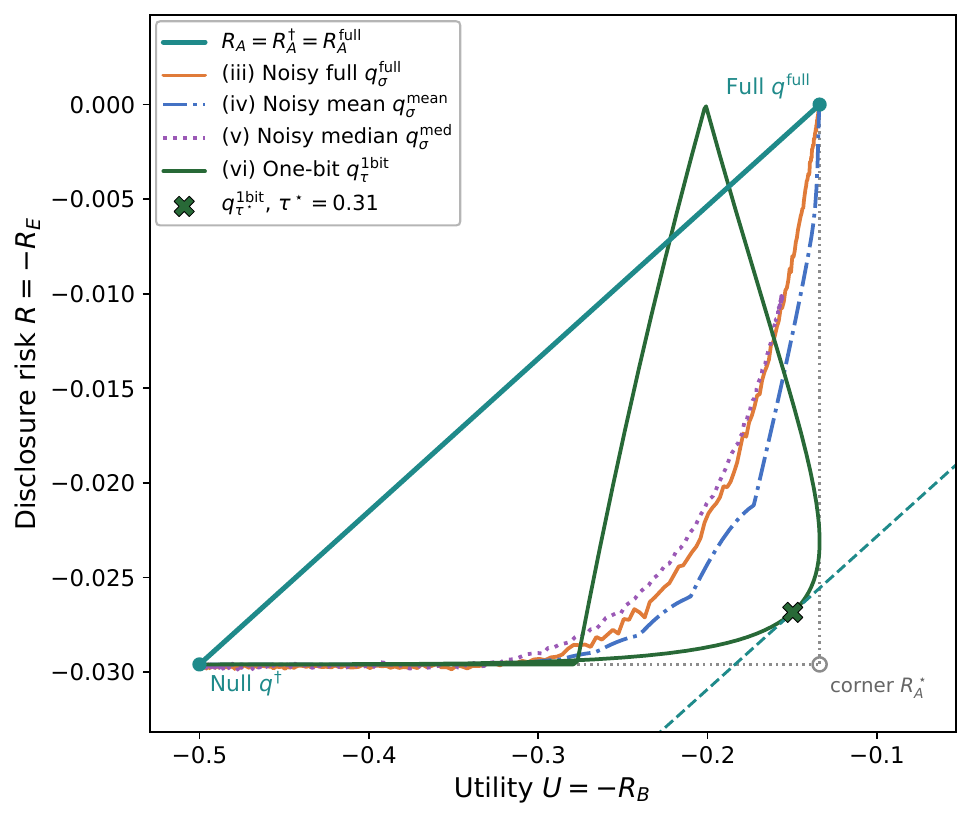}
\caption{Risk--utility map, multiple-Eve criterion. Here the disclosure-risk axis spans a small range ($R_E\lesssim0.03$) and the level line of $R_A$ is steep ($\lambda\approx12.4$).}
\label{fig:supp-rumap-neither}
\end{figure}

We reproduce here, using the numerical results of Example~2, the risk--utility reading introduced in the main text: each mechanism is plotted in the plane of utility $U=-R_B$ (horizontal) and
disclosure risk $R=-R_E$ (vertical). The full release sits at the top right, the null release at the
bottom left, and the calibrated level line of $R_A$ (slope $1/\lambda$) passes through both; $R_A$ decreases towards the bottom right, so mechanisms whose curves lie \emph{below} the line improve on the trivial baseline. Each map superimposes all mechanisms, the additive-noise releases traced as $\sigma$ varies and the one-bit release as $\tau$ varies. In addition to the calibrated line we draw the parallel level line of $R_A$ through the \emph{optimal} operating point $q^\star$ (the mechanism and parameter minimising $R_A$, marked by a star), whose intercept gives the best attainable $R_A$. A dotted vertical line dropped from Full (constant utility $U=-R_B^{\mathrm{full}}$) and a dotted horizontal line drawn from Null (constant disclosure risk $R=-R_E^{\dagger}$) meet at the \emph{theoretical corner} (circle), the point $(U_{\mathrm{full}},R_{\mathrm{null}})$ attaining $R_A^\star=2R_B^{\mathrm{full}}-\min(\pi,1-\pi)$: keeping Full's utility while lowering the risk to Null's level. The optimum $q^\star$ coincides with the corner \emph{exactly when it is reachable} (the maximum adversary; see the risk--utility maps in the main text) and falls short of it otherwise (the mean adversary and the multiple-Eve criterion). For brevity only the multiple-Eve map is reproduced here (Figure~\ref{fig:supp-rumap-neither}); the mean and maximum maps are shown in the main text. Parameters are $\sigma_0=1$, $c_{\mathrm{mean}}=\sigma_0/2$, and an outlier threshold $c_{\max}=2\sigma_0$.

\paragraph{Reproducibility.}
All numerical implementations are written in Python using JAX just-in-time compilation. The full source code used to generate the results and figures in this paper is publicly available at
\url{https://github.com/antoineluciano/BayesianAdversarialPrivacy}.